\documentclass[12pt]{article}

\usepackage{natbib}
\usepackage{amssymb,latexsym}
\usepackage{graphicx,amsmath,amsfonts,amsthm}
\usepackage{color}
\usepackage{graphicx}

\begin{document}
\newtheorem{lemma}{Lemma}[section]
\newtheorem{proposition}{Proposition}[section]
\newtheorem{observation}{Observation}
\newtheorem{theorem}{Theorem}
\newtheorem{definition}{Definition}
\newtheorem{corollary}{Corollary}[section]
\newtheorem{fact}{Fact}[section]
\renewcommand{\theequation}{\thesection.\arabic{equation}}
\newcommand{\cp}{{\cal P}}
\newcommand{\ca}{{\cal A}}
\newcommand{\cb}{{\cal B}}
\newcommand{\cC}{{\cal C}}
\newcommand{\cd}{{\cal D}}
\newcommand{\cR}{{\cal R}}
\newcommand{\cg}{{\cal G}}
\newcommand{\cs}{{\cal S}}
\newcommand{\cl}{{\cal L}}
\newcommand{\cc}{{\cal C}}
\newcommand{\ch}{{\cal H}}
\newcommand{\cv}{{\cal V}}
\newcommand{\cvd}{{\cal V}^{\cd}}
\newcommand{\cvdxy}{{\cal V}^{\cd}_{xy}}
\newcommand{\cu}{{\cal U}}
\newcommand{\cw}{{\cal W}}
\newcommand{\cls}{{\cal LS}}
\newcommand{\Lin}{{\mathcal Lin}(A)}
\newcommand{\hs}{\hspace{1mm}}
\newcommand{\dN}{{\mathbb{N}}}
\newcommand{\Natur}{\dN}
\newcommand{\dfn}{\sf\em}
\def\row#1#2{{#1}_1,\ldots ,{#1}_{#2}}
\thispagestyle{empty}
\vspace*{12mm}
\renewcommand{\thefootnote}{*}
\begin{center}
{\LARGE{\bf Condorcet Domains, Median Graphs and \vspace{3mm}\\
the Single-Crossing Property}} {\large \footnote{We are grateful
to Dominik Peters for very useful comments on the
first version of this paper. Arkadii Slinko was supported by the Marsden Fund grant UOA 1420. {\em More acknowledgements to be added.}}}
\vspace{3mm}
\end{center}
\renewcommand{\thefootnote}{\arabic{footnote}}
\setcounter{footnote}{0}
\vspace{10mm}
\begin{center}
{\sc Clemens Puppe}\vspace{3mm}\\
Department of Economics and Management\vspace{1mm}\\
Karlsruhe Institute of Technology (KIT)\vspace{1mm}\\
D -- 76131 Karlsruhe, Germany\vspace{1mm}\\
clemens.puppe@kit.edu \vspace{5mm}\\
and\vspace{5mm}\\
{\sc Arkadii Slinko}\vspace{3mm}\\
Department of Mathematics\vspace{1mm}\\
The University of Auckland\vspace{1mm}\\
Private Bag 92019\vspace{1mm}\\
Auckland 1142, New Zealand\vspace{1mm}\\
a.slinko@auckland.ac.nz 
\vspace{15mm}\\
\today
%November 2014
\vspace{12mm}
\end{center}
\newpage\noindent
\thispagestyle{empty}
%
%\begin{abstract}
{\bf Abstract.} Condorcet domains are sets of linear orders
with the property that, whenever the preferences of all
voters of a society belong to this set, their majority relation has no cycles. We observe that, without
loss of generality, every such domain can be assumed to be closed in the sense that it
contains the majority relation of every profile with an odd number of voters
whose preferences belong to this domain.

We show that every closed Condorcet domain can be endowed with the structure of a median graph and that, conversely, every
median graph is associated with a closed Condorcet domain (in general, not uniquely).
%The subclass of those 
{\color{black} Condorcet domains that have linear graphs (chains) associated with them are exactly the preference domains with the classical
single-crossing property. As a corollary, we obtain that a domain with the so-called} 
`representative voter property' is either a single-crossing domain or a
very special domains containing exactly four different preference orders whose
associated median graph is a 4-cycle.

Maximality of a Condorcet domain imposes additional restrictions on the associated
median graph. We prove that among all trees only (some) chains can be associated graphs of
maximal Condorcet domains, and we characterize those single-crossing
domains which are maximal Condorcet domains.

Finally, using  the characterization of \cite{NehringPuppe2007b} of
monotone Arrovian aggregation, our analysis yields a rich class
of strategy-proof social choice functions on any 
closed Condorcet domain.\vspace{12mm}\\
{\bf JEL Classification} D71, C72\vspace{1mm}\\
{\bf Keywords:} Social choice, Condorcet domains, acyclic sets of
linear orders, median graphs, single-crossing property, distributive lattice,
Arrovian aggregation, strategy-proof\-ness, intermediate preferences.
%\end{abstract}
%
\setcounter{page}{0}
\newpage

%\vspace*{0mm}
\section{Introduction}
The problem of finding and characterizing preference
domains on which pairwise majority voting never admits cycles---the so-called Condorcet
domains---has a long history in social choice theory. In their seminal
contributions, \cite{Black1948} and \cite{Arrow1951} noticed that the
domain of all linear preference orders that are single-peaked with respect to
some underlying linear spectrum form a Condorcet domain. Later,
\cite{Sen1966} provided a characterization of Condorcet domains in terms of the
well-known condition of value restriction. Since this early work some
progress has been made in understanding the structure
of Condorcet domains; see \cite{Abello91,Abello2004,Cha-Nem1989,DKK:2012,
DanilovK13,PF:1997,PF:2002,GR:2008} for important contributions,
and \cite{mon:survey} for an excellent survey. However, with the exception
of \cite{DanilovK13},
the bulk of the results established in the literature pertains only to the special case of
connected domains. The analysis of \cite{DanilovK13}, on the other
hand, is confined to the case of `normal' Condorcet domains that contain at least
one pair of completely reversed orders.

The present paper provides a unifying general approach
by establishing a close connection between Condorcet domains and
median graphs (see a comprehensive survey about these in \cite{KlavzarMulder})
on the one hand, and the well-studied class
of {\em single-crossing} domains on the other hand; see \cite{Roberts1977,GansSmart1996,Sapor2009},
among others.

First, we observe that if a Condorcet domain of linear preference orders admits a profile with an odd number of voters
whose (transitive) majority relation does not belong to this domain,  then we can add this majority relation to the Condorcet domain
and obtain a larger Condorcet domain.
%
%First, we observe that if one adds to a Condorcet domain of linear preference orders
%the (transitive) majority relation of any profile from this domain with an odd number of voters (if it was not there),
%one obtains again a Condorcet domain. 
We may thus assume without
loss of generality that Condorcet domains are {\em closed} in the sense that
pairwise majority voting among an odd number of individuals always yields
an order within the given domain. In particular, all {\em maximal}
Condorcet domains are necessarily closed. 

The concept of betweenness, introduced for linear orders by  \cite{Kemeny1959}, plays a major role in our analysis.   
An order is {\em between} two orders, or {\em intermediate}, if it agrees
with all binary comparisons in which the two linear orders agree,
(see also~\cite{kem-sne:b:polsci:mathematical-models,Grandmont1978}).
This allows us to associate a graph to any domain of linear orders by calling two linear orders of this domain neighbors if there are no other linear orders in this domain that are between them.

We show that (i) every closed Condorcet domain on a finite set of alternatives equipped with the neighborhood relation
is isomorphic to a median graph, and (ii) for every finite median graph there exists a set
of alternatives and a closed Condorcet domain whose graph is isomorphic to the given median graph.  Importantly,
the median graph corresponding to a Condorcet domain is in general {\em not} a subgraph of the {\em permutohedron}, which is the graph associated
with the universal domain of all strict linear orders (see Section~3 below for detailed explanation)\footnote{The permutohedron was defined in \cite{guilbaud1963analyse} as a polytope. It can also be viewed as a graph in which two vertices are neighbors if they are joined by an edge (which is  the 1-skeleton of that polytope). \cite{mon:survey} calls it {\em permuto\`{e}dre graph}.}.

Our analysis is related to the prior work by \cite{NehringPuppe2007b} and \cite{Demange2012}.
\cite{NehringPuppe2007b} introduce a general notion of  a median space and demonstrate 
its usefulness in aggregation theory; we comment on the relation to their work in detail below.
\cite{Demange2012} shows that, if it is possible to assign linear orders of a certain domain to
the vertices of a median graph in such a way that intermediate orders lie on a shortest path
in the graph, then the majority relation of any profile with preferences from this domain admits
no cycles. However, \cite{Demange2012} takes both the median
graph and the preference profile as given and does not address the issue under which
conditions the required construction can indeed be carried out. Part (i) of our aforementioned result
shows that, in fact, {\em every} closed Condorcet domain (with its respective neighborhood
relation) is isomorphic to a median graph. The analysis of the present paper thus
demonstrates that Demange's requirement of the  existence of a median graph is not only sufficient but also
necessary for a domain of preferences to admit an acyclic majority relation for all profiles. % Importantly,
%this median graph is in general {\em not} a subgraph of the {\em permutohedron}, which is the graph corresponding
%to the natural neighborhood relation on the universal domain of all strict orders (see
%Section~3 below for detailed explanation).

The close connection between Condorcet domains and median graphs established here 
%is useful in particular because it 
allows one to apply results from
the theory of median graphs to shed further light on the structure of
Condorcet domains. For instance, as a simple corollary of our analysis we obtain
that all closed Condorcet domains that contain two completely reversed linear orders
are {\em distributive lattices}. This important fact was first established for
connected maximal Condorcet domains
by \cite{Cha-Nem1989} and \cite{Abello91}, and then for all maximal
Condorcet domains by \cite{DanilovK13} who used a different technique. Here we show this for the
much larger class of all closed (but not necessarily maximal) Condorcet domains. 
%that contain two completely reversed orders.

While all median graphs give rise to closed Condorcet domains, 
%it is clearly not true that every median graph corresponds 
not all of them  correspond to a {\em maximal} Condorcet domain. It turns out
that in fact certain {\em types} of median graphs never enable maximality of the respective
Condorcet domains. In particular,
we prove that among all trees only (some) chains are associated with maximal Condorcet
domains.

Condorcet domains, whose median graphs are chains, have been studied
quite extensively in economics under the name of single-crossing domains. %Using
%the isomorphism between the shortest path (geodesic) betweenness of the
%underlying median graph and the notion of intermediateness of preferences, 
They are characterized by  the {\em single-crossing
property which} stipulates that the orders of the domain can be arranged in a chain so that, for any
ordered pair of alternatives, the set of
all orders that rank one alternative strictly above the other form an interval in this chain.
It is well-known that single-crossing domains have the {\em representative voter
property} (cf.~\cite{Roth91}), i.e.,~in any profile with an odd number of voters whose
preferences belong to the given
domain there is one voter whose preference order coincides with the
majority relation. We show here that the representative voter
property essentially characterizes the class of single-crossing domains.\footnote{The only exceptions
are the domains with exactly four elements such that the associated graph is a $4$-cycle;
these also satisfy the representative voter property.}

%The maximal chains correspond to those single-crossing domains that are not contained
%in any strictly larger single-crossing domain (in particular, these must then contain two
%completely reversed orders). 
A maximal  single-crossing domain must obviously contain two completely reversed orders. 
Interestingly, not all maximal single-crossing domains are maximal Condorcet domains,
i.e.,~typically it is possible to add further preference orders to a maximal single-crossing domain 
 without generating cycles in the majority relation of any profile.
Here, we provide a simple necessary and sufficient condition of when a maximal
single-crossing domain is also a maximal Condorcet domain. The condition requires that
the `switched' pairs of alternatives associated with any two consecutive orders of the domain
have one element in common.

Our analysis immediately implies that the closed Condorcet domains are exactly the
{\em median stable} subsets of the space of all linear orders on a given set of
alternatives. It thus follows from the results of \cite{NehringPuppe2007b,NehringP2010}
that closed Condorcet domains not only enable pairwise majority voting
as a consistent aggregation method but in fact admit a wide range of further
aggregation rules satisfying Arrow's independence condition. We adapt their
characterization of all monotone Arrovian aggregators to the case of Condorcet
domains  of linear orders and
% observe its close connection to the well-known fact that the family
%of all convex subsets of vertices of a median graph possesses the {\em Helly property}, which states %that any
%subfamily of pairwise intersecting convex subsets has a non-empty %common 
%intersection.
%\footnote{{\color{red} In fact, \cite{Bandelt89} showed that the Helly property for
%convex sets together with the absence of subgraphs isomorphic to $K_{2,3}$
%characterize median graphs.}}.
%Finally, we
show that every monotone Arrovian aggregator on a closed Condorcet
domain induces a strategy-proof social choice function on the same domain. For each
Condorcet domain, we thus obtain a rich class of strategy-proof voting rules,
of which the rules identified by \cite{Sapor2009} for single-crossing domains are
special cases. 

The remainder of the paper is organized as follows. In the following Section~2, we
introduce the concept of a Condorcet domain and observe some of its fundamental properties.
In particular, we show that closed Condorcet domains are exactly the median
stable subsets of the space of all linear preference orders on a given set of
alternatives. Section~3
introduces median graphs and states our main result establishing
the correspondence between closed Condorcet domains and median graphs.
Section~4 provides the characterization of single-crossing domains, and
discusses a weaker version of the single-crossing property, namely single-crossingness
on trees. Section~5
addresses maximality of Condorcet domains, an issue that has already received
attention in the literature (cf.~\cite{mon:survey}). In particular, we prove
that trees different from chains are never associated with maximal Condorcet domains,
and we characterize the chains (i.e.~single-crossing domains) that correspond to
maximal Condorcet domains.
Any maximal single-crossing domain necessarily contains two completely reversed orders;
this does not hold for all maximal Condorcet domains as we show by means of an example.
In Section~6, we adapt the characterization of all monotone Arrovian aggregators
obtained in \cite{NehringPuppe2007b} to the case of (closed) Condorcet domains
and show how it induces a large class of strategy-proof social choice functions
on any such domain. An appendix briefly reviews the theory of geometric
interval operators (see, e.g.~\cite{Vel}) and contains the proof of the central Lemma \ref{equivalence} using the so-called triangle condition introduced by \cite{BandeltChepoi1996}.

%%%%%%%%%%%%%%%%%%%%%%%%%%%%%%%%%%%%%%%%%%%%%%
%%%%%%%%%%%%%%%%%%%%%%%%%%%%%%%%%%%%%%%%%%%%%%
%%%%%%%%%%%%%%%%%%%%%%%%%%%%%%%%%%%%%%%%%%%%%%
%%%      SECTION 2
%%%%%%%%%%%%%%%%%%%%%%%%%%%%%%%%%%%%%%%%%%%%%%
%%%%%%%%%%%%%%%%%%%%%%%%%%%%%%%%%%%%%%%%%%%%%%
%%%%%%%%%%%%%%%%%%%%%%%%%%%%%%%%%%%%%%%%%%%%%%

\section{Closed Condorcet Domains as the Median Stable \\ Subsets of the Universal Domain}

In this section, we show that the class of domains that are closed under taking the majority
relation for all profiles with an odd number of voters are precisely the median stable subsets of
the space of linear orders endowed with the Kemeny betweenness relation. {\color{black} A majority of results in this section are not original. In Theorem 1 we gathered a number of classical characterisations of Condorcet domains and Theorem 2 can be derived from  the analysis in \cite{NehringPuppe2007b} who investigated a more general class of median spaces.}
%To the best of our knowledge, the central notion of this paper, {\color{black} that of} a
%closed Condorcet domain, has not yet been formally introduced in the literature. On
%the other hand, the main results of this section are {\color{black} straightforward adaptions
%of the analysis in \cite{NehringPuppe2007b} 
%(see} the note at the end of this section). 
Nevertheless, we believe that our exposition and the short proof of the fundamental characterization result, Theorem \ref{cc-iff-med} below, will help to clarify and unify several
different approaches in the literature.

\subsection{Condorcet Domains}
Consider a finite set of alternatives $X$ and the set $\cR (X)$ of all {\dfn (strict) linear 
orders} (i.e.,~complete, transitive and antisymmetric binary relations) on $X$ which we will refer to as the {\dfn universal domain}.
Any subset $\cd\subseteq \cR(X)$ will be called a {\dfn domain of preferences} or simply a {\dfn domain}. A
{\dfn profile $\rho=(\row Rn)$ over $\cd$} is an element of the Cartesian product $\cd^n$ for some number $n\in\Natur$
of `voters', where the linear order $R_i$ represents the
preferences of the $i$th voter over the alternatives from $X$.
A profile with an odd number of voters will simply be referred to as an {\dfn odd profile}.
Frequently, we will denote linear orders simply by listing the alternatives in the order of decreasing preference, e.g., a linear
order that ranks $a$ first, $b$ second, $c$ third, etc., is denoted by $abc\ldots$.

The {\dfn majority relation} associated with a profile $\rho$ is the binary relation
$P_{\rho}^\text{maj}$ on $X$ such that $x  P_{\rho}^\text{maj}  y$ if and only if
more than half of the voters rank $x$ above $y$. Note that, according to this definition, the
majority relation is asymmetric and, for any odd profile $\rho$ and any two distinct
alternatives $x,y\in X$, we have either $x  P_{\rho}^\text{maj}  y$ or $y P_{\rho}^\text{maj}  x$.
%
\iffalse
Throughout, we will use the following notational conventions.
Since linear orders are assumed to be reflexive, they are generically denoted by the
capital letter `$R$,' and their asymmetric part by the capital letter `$P$.' When we
speak of `the' majority relation corresponding to an odd profile $\rho$, we mean either
the asymmetric relation $P_{\rho}^{maj}$ or the corresponding reflexive relation
$R_{\rho}^{maj}:=P_{\rho}^{maj}\cup\{ (x,x)  \mid  x\in X\}$. 
\fi
%
An asymmetric binary relation $P$ is {\dfn acyclic} if there does not exist a subset
$\{ x_1,\ldots,x_m\}\subseteq X$ such that $x_1 P x_2$, $x_2 P x_3$, \ldots , $x_{m-1} P x_m$
and $x_m P x_1$. The class of all domains $\cd\subseteq\cR (X)$ such that, for all $n$, the
majority relation associated with any profile $\rho\in \cd^n$ 
is acyclic has received significant attention in the literature, see the survey of 
\cite{mon:survey} and the references therein. In the following, we will refer to any such
domain as a {\dfn Condorcet domain}.\footnote{\cite{PF:1997} calls them {\em acyclic
sets of linear orders}.}

The following result prepares the ground for our analysis, providing some well-known
characterizations of Condorcet domains (cf.~\cite[p.~142]{mon:survey}). In particular,
condition d) below is Sen's [1966] `value restriction'; condition e) has been introduced
by \cite{Ward1965} as the `absence of a Latin square' (in other terminology, it requires
the absence of a `Condorcet cycle'; cf.~\cite{Condorcet1785} [1785]). 
%%%

\begin{theorem}
\label{prel_facts}
Let $X$ be finite, and let $\cd\subseteq\cR (X)$ be a subset of the
space of all linear orders on $X$. The following statements are equivalent:
\begin{itemize}

\item[{a)}] $\cd$ is a Condorcet domain, i.e.,~the majority relation corresponding to
every profile over $\cd$ is acyclic.

\item[{b)}] For every profile over $\cd$, the corresponding majority relation is a strict partial
order (i.e.,transitive and asymmetric binary relation).

\item[{c)}] For every odd profile over $\cd$, the corresponding majority relation is a linear order,
i.e.,~an element of $\cR (X)$.

\item[{d)}] For every triple $x,y,z \in X$ of pairwise distinct alternatives, there exists one element
in $\{ x,y,z\}$ that is either never ranked first, or never ranked second, or never ranked third 
in all restrictions of the orders in $\cd$ to the set $\{ x,y,z\}$.

\item[{e)}] For no triple $R_1, R_2, R_3\in\cd$, and no triple $x,y,z\in X$ of pairwise
distinct alternatives one has 
$x R_1 y R_1 z$, $y R_2 z R_2 x$ and $z R_3 x R_3 y$ simultaneously.

\end{itemize}
\end{theorem}

%%%

We will say that a Condorcet domain $\cd$ is {\dfn closed} if the majority relation
corresponding to any odd profile over $\cd$ is again an
element of $\cd$, and we will say
that a Condorcet domain $\cd$ is {\dfn maximal} if no Condorcet domain (over the same
set of alternatives) is a proper superset of $\cd$. The following simple observation
will be very useful.

%%%

\begin{lemma} 
\label{closedCdomains}
Let $\cd$ be a Condorcet domain and $R\in\cR (X)$ be
the majority relation corresponding to an odd profile over $\cd$.
Then $\cd\cup\{ R\}$ is again a Condorcet domain. In particular, every Condorcet
domain is contained in a closed Condorcet domain.
\end{lemma}

\noindent {\bf Proof.} By Theorem 1e), it suffices to show that 
$\cd\cup\{ R\}$ does not admit three orders $R_1, R_2, R_3$
and three elements $x,y,z\in X$ such that $xR_1 yR_1 z$, $yR_2 zR_2 x$ and
$zR_3 xR_3 y$. Assume on the contrary that it does; then, evidently, not all three
orders $R_1, R_2, R_3$ belong to $\cd$. Thus, one of them, say $R_3$, is the
majority relation of an odd profile $\rho\in\cd^n$, i.e.~$R_3=R$.
Consider the profile $\rho' =(nR_1, nR_2,\rho) \in\cd^{3n}$ that
consists of $n$ voters having the order $R_1$, $n$ voters having the order $R_2$
and the $n$ voters of the profile $\rho$.
Then the voters of the subprofile $(nR_1, nR_2)$ will
unanimously prefer $y$ to $z$, which
forces the majority relation $P_{\rho' }^\text{maj}$ corresponding to $\rho'$ to
have the same ranking
of $y$ and $z$. At the same time, the voters of this subprofile are evenly
split in the ranking of any
other pair of alternatives from $\{x,y,z\}$. Hence, the majority
relation $P_{\rho' }^\text{maj}$
yields the cycle $z P_{\rho' }^\text{maj} x P_{\rho' }^\text{maj} y P_{\rho' }^\text{maj} z$, in
contradiction to the assumption that $\cd$ is a Condorcet domain.%\vspace{-2mm} \\
\hspace*{1mm} \hfill $\Box$

\medskip

%%%

This observation allows us to concentrate our attention on closed
Condorcet domains without loss of generality, and we do so for the rest of the paper.
Note, in particular, that by
Lemma~\ref{closedCdomains} all maximal Condorcet domains are closed.

\subsection{Betweenness and Median Domains}

The universal domain $\cR (X)$ can be endowed with the  ternary relation of
`Kemeny betweenness' \citep{Kemeny1959}.
According to it, an order $Q$ is {\dfn between} orders $R$ and $R'$ if
$Q\supseteq R\cap R'$, i.e.,~$Q$ agrees with all binary comparisons on which
$R$ and $R'$ agree.\footnote{Some
%(\cite{kem-sne:b:polsci:mathematical-models})
authors such as, e.g., \cite{Grandmont1978} and \cite{Demange2012} refer to orders that
are between two others in this sense as `intermediate' orders.} The set of all
orders that are between $R$ and $R'$ is called the {\dfn interval}
spanned by $R$ and $R'$ and is denoted by $[R,R']$.
This {\dfn interval operator} makes $\cR (X)$ an {\em interval space}
\citep{Vel}.\footnote{ For a closer analysis of this interval operator,
see the appendix.}

A subset $\cd\subseteq\cR (X)$ %of the permutohedron 
is called {\dfn median stable} if, for any triple of elements
$R_1,R_2,R_3\in\cd$, there exists an element $R^\text{med} =
R^\text{med}(R_1,R_2,R_3)\in\cd$, the {\dfn median order} corresponding
to $R_1,R_2,R_3$, such that
\[
R^\text{med}\in  [R_1,R_2]\cap[R_1,R_3]\cap[R_2,R_3].
\]

%%%

\begin{proposition}
\label{uniqueness_prop}
The median order of a triple $R_1, R_2, R_3\in \cR (X)$, if exists, is unique. 
\end{proposition}

\noindent{\bf Proof.} If a triple $R_1, R_2, R_3$ admits two different
median orders, say $R$ and $R'$, these must differ on the ranking of at least one pair of alternatives.
Suppose they disagree on the ranking of $x$ and $y$. In this case, not all three
orders of the triple agree on the ranking of $x$ versus $y$. Hence, exactly two of them,
say $R_1$ and $R_2$, must agree on the
ranking of $x$ versus $y$; but then, either $R$ or $R'$ is not between $R_1$ and $R_2$, a
contradiction.\vspace{-4mm} \\  \hspace*{1mm} \hfill \vspace{2ex} $\Box$

%%%

In the following, we will refer to median stable subsets of $\cR (X)$
as {\dfn median domains}.
Evidently, not every subset of $\cR(X)$ is a median domain; for instance, the universal
domain $\cR(X)$ itself is not a median domain whenever $|X|\geq 3$. This can be verified
by considering any three orders of the form $R_1 = \ldots a \ldots b \ldots c \ldots$,
$R_2 = \ldots b \ldots c \ldots a \ldots$, and $R_3 = \ldots c \ldots a \ldots b \ldots$.
Since any linear order $R$ in $[R_1,R_3]$ has $aRb$,
any linear order $R$ in $[R_1,R_2]$ has $bRc$, and any linear order $R$ in
$[R_2,R_3]$ has $cRa$, we obtain $[R_1,R_2]\cap[R_1,R_3]\cap[R_2,R_3] = \emptyset$
due to the transitivity requirement.
Prominent examples of median domains include
the well-studied single-crossing domains.
\medskip

\noindent {\bf Example 1 (Classical single-crossing domains).} There are several
equivalent descriptions of single-crossing domains (see, e.g.,~\cite{GansSmart1996,Sapor2009}).
The following will be useful
for our purpose. A domain $\cd\subseteq\cR (X)$ is said to have the
{\em single-crossing property}
if $\cd$ can be linearly ordered, say according to $R_1 > R_2 > \ldots > R_m$, so
that, for all pairs $x,y$ of distinct elements of $X$, the sets
$\{ R_j\in\cd \mid xR_j y\}$ and $\{ R_j\in\cd \mid yR_j x\}$ are connected in the
order $>$. Thus, for each pair $x,y$ of distinct elements, there is exactly one `cut-off' order
$R_k$ such that either (i) $xR_j y$ for all $j\leq k$ and $yR_jx$ for all $j>k$, or (ii) $yR_j x$
for all $j\leq k$ and $xR_jy$ for all $j>k$. 
It is easily verified that, for any triple with $R_i>R_j>R_k$ the median order
exists and coincides with the `middle' order, i.e., $R^\text{med}(R_i,R_j,R_k) =R_j$.
\medskip

The close connection between Condorcet domains and median domains to be
established in Theorem 2
below stems from the following simple but fundamental observation
(cf.~\cite[Cor.~5]{NehringPuppe2007b}).
\medskip

\noindent{\bf Observation.} {\em A triple $R_1, R_2, R_3\in \cR (X)$ admits a
median order if and only if the majority relation of the profile $\rho =(R_1, R_2, R_3)$ is acyclic,
in which case the median
order $R^\text{med}(R_1,R_2,R_3)$ and the majority relation of $\rho$ coincide.} \vspace{1ex}

\noindent {\bf Proof.} If the majority relation $P_{\rho }^\text{maj}$ is acyclic, and hence is an element
of $\cR (X)$, it belongs to each interval $[R_i,R_j]$ for all distinct $i,j\in\{ 1,2,3\}$.
Indeed, if  both $R_i$ and $R_j$ rank $x$ higher than $y$, then so does the majority relation.
Conversely, if $R$ is the median of the triple
$R_1, R_2, R_3$, then for any pair $x,y\in X$, at least two orders from this triple
agree on ranking of $x$ and $y$. Then $R$ must agree with them,
hence it is the majority relation for the  profile
$\rho =(R_1, R_2, R_3)$.\vspace{-4mm} \\ 
\hspace*{1mm} \hfill $\Box$
%%%
\begin{corollary}  
\label{cc-is-med-corollary}
Any  closed Condorcet domain is a median domain.
\end{corollary}

\noindent {\bf Proof.} Suppose $\cd$ is a closed Condorcet domain, and let
$R_1, R_2, R_3$ be any triple of orders from $\cd$. The majority relation $R$
corresponding to the profile $(R_1, R_2, R_3)\in\cd^3$ by Theorem~1c) is an
element of $\cR (X)$, 
and by the assumed closedness it is in fact an element of $\cd$.
By the  preceding observation, $R$ is the median order of the triple $R_1, R_2, R_3$.
\vspace{-4mm} \\ \hspace*{1mm} \hfill $\Box$
%\pagebreak
%%%
\medskip

%%%

A subset $\cc\subseteq\cd$ of a domain $\cd\subseteq\cR (X)$ will be called {\dfn convex}
if $\cc$ contains with any pair $R,R'\in\cc$ the entire interval spanned by $R$ and $R'$, that
is, $\cc$ is convex if
\[
\{ R,R'\}\subseteq \cc \  \Rightarrow  \ [ R, R'] \subseteq \cc.
\]

A family ${\mathbb F}$ of subsets of a set is said to have the {\dfn Helly property} if
the sets in any subfamily ${\mathbb F}'\subseteq {\mathbb F}$ have a non-empty intersection
whenever their pairwise  intersections are non-empty, i.e.,
if $\cc\cap\cc'\neq\emptyset$ for each pair $\cc,\cc'\in{\mathbb F}'$ implies
$\cap \, {\mathbb F}'\neq\emptyset$. For us this property will be important
when ${\mathbb F}$ is the set of all convex subsets of a domain $\cd\subseteq\cR (X)$.

%%%

\begin{proposition}[\bf Helly property and median domains]
\label{helly-prop}
A domain $\cd$ is a median domain if and only if 
$\cd$ has the Helly property for convex subsets.
\end{proposition}

\noindent {\bf Proof.} Let $\cd$ be median domain and ${\mathbb F}$ be a family of convex
subsets with pairwise non-empty intersection. We proceed by induction over $m=|{\mathbb F}|$.
If $m=2$, there is nothing to prove, thus let $m=3$, i.e., ${\mathbb F} = \{ \cc_1,\cc_2,\cc_3\}$.
Choose any triple of orders $R_1\in\cc_1\cap\cc_2$,
$R_2\in\cc_2\cap\cc_3$ and $R_3\in\cc_3\cap\cc_1$, and consider the median order
$R = R^\text{med}(R_1,R_2,R_3)$. By convexity of the sets $\cc_1,\cc_2,\cc_3$ we
have $R\in \cc_1\cap\cc_2\cap\cc_3$ which, in particular, shows that $\cap\, {\mathbb F}$
is non-empty.
Now consider ${\mathbb F} = \{ \cc_1,\ldots,\cc_m\}$ with $m>3$ elements, and assume that
the assertion holds for all families with less than $m$ elements.
Then, the family $\{ \cc_1, \cc_2, \cc_3\cap \ldots \cap\cc_m\}$ constitutes a family of
three convex subsets with pairwise non-empty intersections. By the preceding argument,
we thus have $\cap\, {\mathbb F}\neq\emptyset$.

Conversely, consider a domain $\cd$ such that any family of convex subsets of $\cd$
has the Helly property. Consider any three orders $R_1,R_2,R_3\in\cd$. Since, evidently,
all intervals are convex, the Helly property applied to the intervals $[R_1,R_2]$, $[R_1,R_3]$,
$[R_2,R_3]$ implies the existence of a median.\vspace{-4mm} \\  
\hspace*{1mm} \hfill \vspace{1ex} $\Box$

%%%

For any domain $\cd$ and any pair $x,y\in X$ of alternatives, denote by $\cvdxy$
the set of orders in $\cd$ that rank $x$ above $y$, i.e.,
\[
\cvdxy := \{ R\in\cd \mid xRy\}.
\]
Note that, for all distinct $x,y\in X$, the sets $\cvdxy$ and $\cvd_{yx}$ form a
partition of $\cd$. Also observe that the sets of the form $\cvd_{xy}$ are convex
for all pairs $x,y\in X$.
We will now use the Helly property applied to this family of convex sets
to show that every median domain is a closed Condorcet domain.
The following is the main result of this section.

\begin{theorem}  
\label{cc-iff-med}
The classes of median domains and closed Condorcet domains coincide, i.e., a domain
is a median domain if and only if it is a closed Condorcet domain.
\end{theorem}  

\noindent {\bf Proof.} In the light of Corollary~\ref{cc-is-med-corollary}, it suffices to show
that every median domain is a closed Condorcet domain. 
Thus, let $\cd$ be a median domain and consider an odd profile
$\rho =(\row Rn)\in\cd^n$. For any two alternatives $x,y\in X$, let
$\cu_{xy}=\{R_i \mid xR_iy \}$, and observe that obviously, $\cu_{xy}\subseteq\cvd_{xy}$.
Let $z,w$ also be alternatives in $X$, not necessarily distinct from $x$ and $y$. If
$xP^\text{maj}_{\rho}y$ and $zP^\text{maj}_{\rho}w$, then
$\cu_{xy}\cap\, \cu_{zw}\neq\emptyset$ and hence
$\cvd_{xy}\cap \cvd_{zw}\neq\emptyset$.
By Proposition~\ref{helly-prop} we have
\[
\bigcap_{x,y\in X \, : \, xP^\text{maj}_{\rho}y}  \cvd_{xy}\ne \emptyset,
\]
hence there is a linear order in $\cd$ which coincides with the majority
relation of $\rho$. \vspace{-4mm} \\ 
\hspace*{1mm}  \hfill \vspace{2ex} $\Box$

%%%

%\noindent
%{\bf Relation to the literature.} Median domains as studied here
%belong to the general class of median spaces analyzed in \cite{NehringPuppe2007b}.
%Theorem \ref{cc-iff-med} can be derived from their analysis, and in fact they
%use the Helly property without referring to it by this name.

%%%%%%%%%%%%%%%%%%%%%%%%%%%%%%%%%%%%%%%%%%%%%%
%%%%%%%%%%%%%%%%%%%%%%%%%%%%%%%%%%%%%%%%%%%%%%
%%%%%%%%%%%%%%%%%%%%%%%%%%%%%%%%%%%%%%%%%%%%%%
%%%      SECTION 3
%%%%%%%%%%%%%%%%%%%%%%%%%%%%%%%%%%%%%%%%%%%%%%
%%%%%%%%%%%%%%%%%%%%%%%%%%%%%%%%%%%%%%%%%%%%%%
%%%%%%%%%%%%%%%%%%%%%%%%%%%%%%%%%%%%%%%%%%%%%%

\section{Closed Condorcet Domains and Median Graphs}
\setcounter{equation}{0}

The close relation between the majority relation and the median operator in closed Condorcet domains suggests
that there exists a deeper structural connection between closed Condorcet domains and
median graphs.\footnote{The term `median graph' was coined by \cite{Nebesky1971}; for a comprehensive
survey on median graphs see \cite{KlavzarMulder}.} The details of this connection are worked out
in this section. We start in
Subsection~3.1 with some basic facts about median graphs. In Subsection~3.2 we
prove that the graph associated with any closed Condorcet domain is a median graph. In Subsection~3.3 we show that, conversely, for every median graph one can construct a (non-unique) closed Condorcet domain whose associated graph is isomorphic to the given graph.

\subsection{Median Graphs}

Let $\Gamma =(V,E)$ be a connected graph. The {\dfn distance} $d(u,v)$ between two
vertices $u,v\in V$ is the smallest number of edges that a path connecting $u$
and $v$ may contain. While the distance is uniquely
defined, there may be several shortest paths from $u$ to $v$. We say that a
vertex $w$ is {\dfn geodesically between} the vertices $u$ and $v$ if $w$ lies on a
shortest path that connects
$u$ and $v$ or, equivalently, if $d(u,v)=d(u,w)+d(w,v)$. 
A {\dfn (geodesically) convex} set in a graph $\Gamma =(V,E)$ is a subset
$C\subseteq V$ such that for any two
vertices $u,v\in C$ all vertices on every shortest path between $u$ and $v$ in
$\Gamma$ lie in $C$. A connected graph $\Gamma =(V,E)$ is called a
{\dfn median graph} if, for any three vertices $u,v,w\in V$,
there is a \emph{unique} vertex $\text{med} (u,v,w)\in V$ which lies simultaneously on some
shortest paths from $u$ to $v$, from $u$ to $w$ and from $v$ to $w$.

To characterize the structure of an arbitrary median graph we recall the concept of
{\em convex expansion}. For any two subsets $S,T\subseteq V$ of the set of vertices of
the graph $\Gamma$, let $E(S,T)\subseteq E$ denote the
set of edges that connect vertices in $S$ and vertices in $T$.

\begin{definition}
Let $\Gamma =(V,E)$ be a graph. 
Let $W_1,W_2\subset V$ be two subsets with a non-empty intersection
$W_1\cap W_2\ne \emptyset$ such that $W_1\cup W_2=V$ and
$E(W_1\setminus W_2,W_2\setminus W_1)=\emptyset$.
The {\dfn expansion} of $\Gamma$ with respect to $W_1$ and $W_2$ is the
graph $\Gamma'$ constructed as follows:
\begin{itemize}
\item  each vertex $v\in W_1\cap W_2$ is replaced by two vertices $v^1$, $v^2$
joined by an edge;
\item  $v^1$ is joined to all the neighbors of $v$ in $W_1\setminus W_2$ and $v^2$ is
joined to all the neighbors of $v$ in $W_2\setminus W_1$;
\item if $v,w\in W_1\cap W_2$ and $vw\in E$, then $v^1$ is joined to $w^1$ and
$v^2$ is joined to $w^2$;
\item if $v,w\in W_1\setminus W_2$ or if $v,w\in W_2\setminus W_1$, they will be
joined by an edge in $\Gamma'$ if and only if they were joined in $\Gamma$;
if $v\in W_1\setminus W_2$ and $w\in W_2\setminus W_1$, they remain not joined in $\Gamma'$.
\end{itemize}
If $W_1$ and $W_2$ are convex, then $\Gamma'$ will be called
a {\dfn convex expansion} of $\Gamma$. 
\end{definition}

\medskip

\begin{figure}
\centering
\includegraphics[width=14cm]{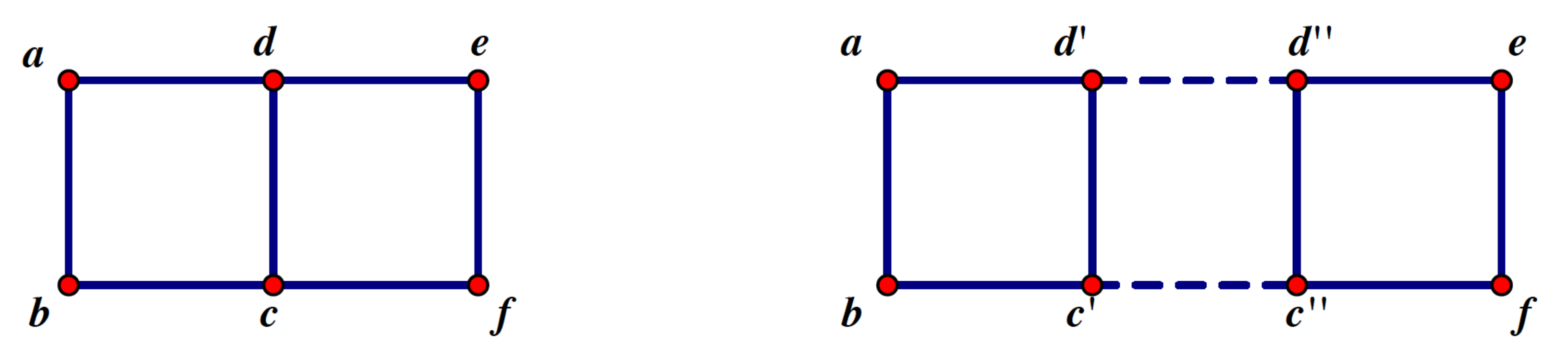}
\caption{Convex expansion of a median graph}
\label{Fig1}
\end{figure}

\noindent {\bf Example 2 (Convex expansion).}
In the graph $\Gamma$ shown on the left of Figure~\ref{Fig1} we set
$W_1=\{a,b,c,d\}$ and $W_2=\{c,d,e,f\}$. These are convex and their
intersection $W_1\cap W_2=\{c,d\}$ is not empty. On the right we see the
graph $\Gamma'$ obtained by the convex expansion of $\Gamma$ with
respect to $W_1$ and $W_2$.

%\vspace{4mm}

%\begin{center} 
%\includegraphics[width=14cm]{Fig-Expansion}
%\end{center}
%\begin{center}
%{\em Figure 1: Convex expansion of a median graph}
%\vspace{3mm}
%\end{center}

The following important theorem about median graphs is due to \cite{mulder}.

\begin{theorem}[\bf Mulder's convex expansion theorem]
A graph is median if and only if it can be obtained from a trivial one-vertex
graph by repeated convex expansions. 
\end{theorem}

%%%%%%%%%%%%%%%%%%%%%%%%%%%%%%%%%%%%%%%%%%%%%
%   SUBSECTION 3.2
%%%%%%%%%%%%%%%%%%%%%%%%%%%%%%%%%%%%%%%%%%%%%

\subsection{Every Closed Condorcet Domain Induces a Median Graph}
\label{bgbd}

As we already said in the introduction, with every domain $\cd\subseteq\cR (X)$ one can
associate a graph
$\Gamma_{\cd}$ on $\cd$ as follows. Two distinct orders
$R,R'\in\cd$ are said to be {\dfn neighbors in $\cd$}, or simply {\dfn $\cd$-neighbors},
if $[R,R']\cap\cd=\{ R,R'\}$.
Define $\Gamma_{\cd}$ to be the (undirected) graph on $\cd$ that connects each pair of 
$\cd$-neighbors by an edge. We say that $\Gamma_{\cd}$ is the {\dfn  associated graph} of $\cd$.
We recall that the graph associated with the universal domain $\cR(X)$ is called
the {\dfn permutohedron}.
Note that the graph $\Gamma_{\cd}$ is always connected, i.e., any two orders in
$\cd$ are connected by a path in $\Gamma_{\cd}$. Moreover,
any two $\cR (X)$-neighbors $R,R'$ are always $\cd$-neighbors whenever $R,R'\in\cd$.
However, two $\cd$-neighbors need not be
$\cR (X)$-neighbors, so, if $\cd\ne \cR(X)$, the associated graph $\Gamma_{\cd}$ need not be a subgraph of the permutohedron. 
If it is, the domain $\cd$ is called {\dfn connected}.\footnote{{\color{black} Since connectedness
of a domain is thus not the same as connectedness of the associated graph, the terminology
might be a bit confusing. But it is widely accepted in the literature \citep{mon:survey}.}}

We will now define another concept of betweenness, different from the Kemeny one. For any $\cd\subseteq\cR (X)$, the order $Q\in\cd$ is
{\dfn $\Gamma_{\cd}$-geodesically
between} the orders $R, R'\in\cd$ if $Q$ lies on a shortest
$\Gamma_{\cd}$-path that connects
$R$ and $R'$. %Finally, say that a domain $\cd$ itself is {\dfn connected} if the graph
%$\Gamma_{\cd}$ is a subgraph
%of $\Gamma_{\cR (X)}$, i.e.,~if all $\cd$-neighbors differ in the ranking of exactly {\em one}
%pair of alternatives.

\begin{figure}[h]
\centering
\includegraphics[width=9cm]{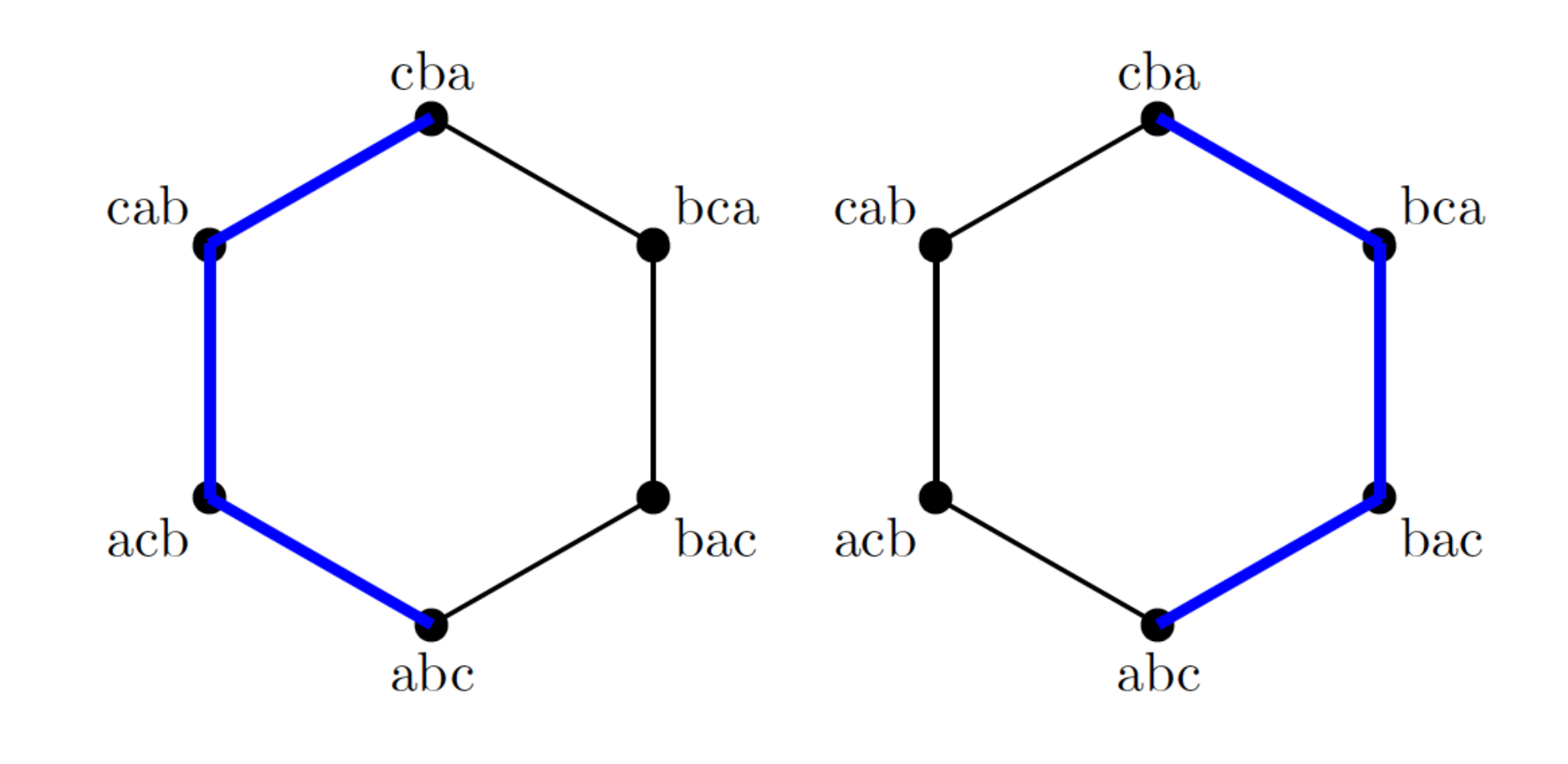}
\caption{Two single-crossing domains on $\{ a,b,c\}$ with their associated graphs}
\label{Fig2}
\end{figure}

As a first example, consider Figure~\ref{Fig2} which depicts two single-crossing domains on the
set $X=\{ a,b,c\}$ with their corresponding graphs.
The domain on the left is $\cd_1=\{ abc, acb, cab, cba\}$, the domain
on the right is $\cd_2=\{ abc, bac, bca, cba\}$. 

%\vspace{0mm}
%\begin{center} 
%\hspace*{-8mm}\includegraphics[width=9cm]{Fig-TwoSingleCross} \vspace{0mm} \\
%
%{\em Figure 2: Two single-crossing domains on $\{ a,b,c\}$ with their associated graphs}
%\vspace{2mm}
%\end{center}
%

Note that the permutohedron is given by the
entire $6$-cycle, and that the graphs associated with the two single-crossing domains
are linear subgraphs of the permutohedron. In particular, both domains are connected and
the Kemeny betweenness relation on $\cd_1$ and $\cd_2$ %defined in Section~2.2 above
translates into the geodesic betweenness of their associated graphs.\par\smallskip

\begin{figure}[h]
\centering
\includegraphics[width=12cm]{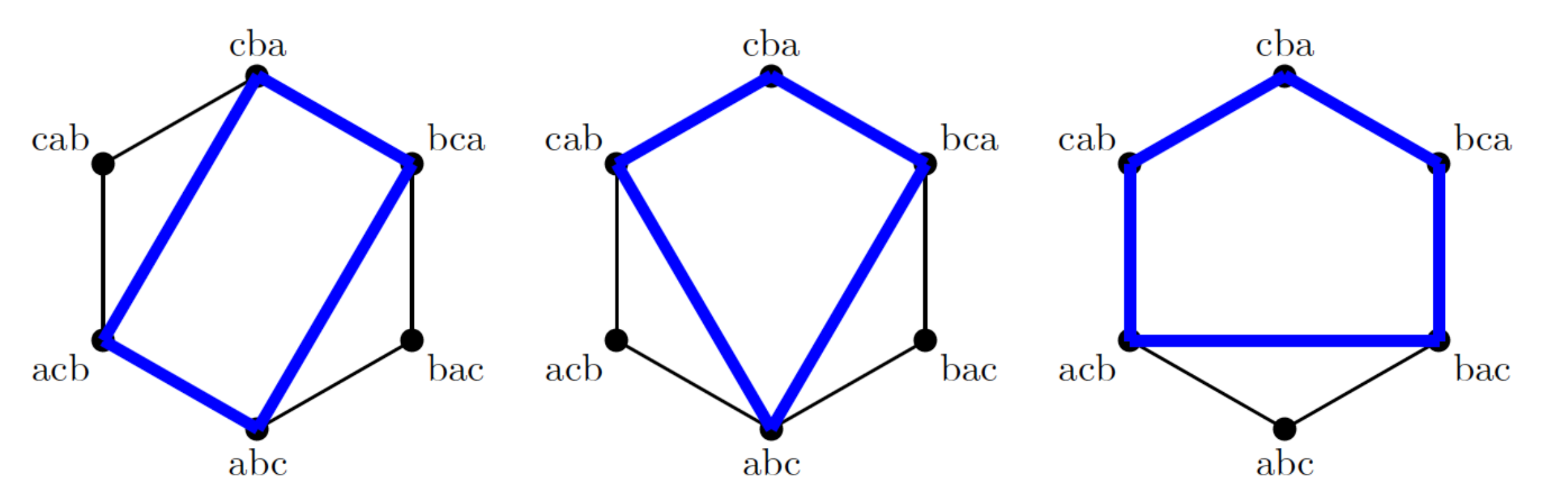}
\caption{Three domains on $\{ a,b,c\}$ with their associated graphs}
\label{Fig3}
\end{figure}

It is important to note that %these features are specific to the particular domains depicted in Fig.~2 and do not hold generally.
for an arbitrary domain $\cd$ both of these properties may not be true:
neither is $\Gamma_{\cd}$ in general a subgraph of $\Gamma_{\cR (X)}$, nor do
the Kemeny betweenness on $\cd$ and the geodesic betweenness on $\Gamma_{\cd}$
correspond to each other. 
To illustrate this, consider the three
domains with their associated graphs
depicted in Figure~\ref{Fig3}. Evidently, none of the three graphs is a subgraph of the permutohedron,
hence none of the corresponding domains is connected. As is easily verified, the domain
$\{ abc, acb, cba, bca\}$ on the left of Figure~\ref{Fig3} is a median domain, but the two other
domains are not; for instance, the domain $\{ abc, cab, cba, bca\}$ in the middle of Figure~\ref{Fig3}
contains a cyclic triple of orders $abc$, $bca$ and $cab$, and the $5$-element
domain $\{ acb, cab, cba, bca, bac\}$ on the right of Figure~\ref{Fig3} contains the cyclic triple
$acb$, $cba$ and $bac$.

%%%
%\begin{center} 
%\includegraphics[width=12cm]{Fig-ThreeDomains} \vspace{2mm} \\
%
%{\em Figure 3: Three domains on $\{ a,b,c\}$ with their corresponding graphs}
%\vspace{2mm}
%\end{center}

%
In case of the domain on the left-hand side, the Kemeny betweenness on $\cd$
and the geodesic betweenness on the associated graph $\Gamma_{\cd}$ agree.
By contrast, for the domain in the middle we have $abc \not\in [cab,bca]$ but
evidently $abc$ is geodesically between the orders $cab$ and $bca$
in the associated graph. On the other hand, in this example the Kemeny betweenness implies
geodesic betweenness. But this also does not hold in general, as the domain on the
right of Figure~\ref{Fig3} shows: here, both orders $cab$ and $cba$ are elements of $[acb, bca]$,
but neither of them is geodesically between $acb$ and $bca$.

%%%
\medskip

 For any domain $\cd\subseteq\cR (X)$, we now have two
betweenness relations on $\cd$:
 the Kemeny betweenness and the geodesic betweenness in the associated
graph $\Gamma_\cd$. We will show that in important cases they coincide, in which case
for all $R,R',Q\in\cd$, we have
\begin{equation}
\label{equivalence_equation}
Q\in [R,R'] \ \Leftrightarrow \ Q \mbox{ is }
\Gamma_{\cd}\mbox{-geodesically between }R\mbox{ and }R'. 
\end{equation}
The following lemma is central to our approach. It states that the Kemeny
betweenness on a domain coincides with the geodesic betweenness of its
associated graph for three natural classes
of domains: (i) median domains, (ii) connected domains, and (iii) 
domains for which the associated graph is acyclic. 

%%%

\begin{lemma} 
\label{equivalence}
The Kemeny betweenness on a domain $\cd\subseteq\cR (X)$ coincides with the geodesic betweenness
of the associated graph $\Gamma_{\cd}$, if one of the following conditions is satisfied:
\begin{itemize}
\item[(i)] $\cd$ is a median domain,
\item[(ii)] $\cd$ is connected,
\item[(iii)] $\Gamma_{\cd}$ is acyclic (i.e., a tree). 
\end{itemize} 
\end{lemma}

The proof of this lemma is provided in the
appendix and invokes a general result on `geometric interval operators'
satisfying the so-called {\em triangle condition} due to \cite{BandeltChepoi1996}. \smallskip
%%%

From Theorem 2 and Lemma~\ref{equivalence}(i), we immediately obtain  the
central result of this section.

%%%
\begin{theorem}
\label{clCo-is-med} 
The associated graph $\Gamma_{\cd}$ of any closed Condorcet domain $\cd$ is a median graph.
Moreover, the Kemeny betweenness relation on %a closed Condorcet domain 
$\cd$ coincides with the geodesic betweenness on $\Gamma_{\cd}$. 
\end{theorem}

We again stress that not all median domains are connected, as exemplified by the domain
on the left-hand-side of Figure~\ref{Fig3}. It is also worth noting that Lemma~\ref{equivalence} does {\em not}
imply that $\cd$ is a median domain whenever the associated graph $\Gamma_{\cd}$
is median graph. A counterexample is the domain $\cd$ in the middle of Figure~\ref{Fig3}, which is
not a median domain despite the fact that its graph $\Gamma_{\cd}$ is a median
graph. However, we have the following corollary.

\begin{corollary} 
Let $\cd\subseteq\cR (X)$ be a connected domain. Then, $\cd$ is a median domain
if and only if the associated graph $\Gamma_{\cd}$ is a median graph.
\end{corollary}
\noindent{\bf Proof.} The associated graph of any median domain is a median graph
by Lemma~\ref{equivalence}(i). Conversely, if $\cd$ is connected the geodesic median
of any triple of vertices with respect to $\Gamma_{\cd}$ is also the median with
respect to the Kemeny betweenness in $\cd$ by Lemma~\ref{equivalence}(ii). Thus, if
$\Gamma_{\cd}$ is a median graph, $\cd$ is a median domain.% \vspace{-4mm} \\ 
\hspace*{1mm} \hfill  $\Box$\medskip

As another important corollary of our analysis, we obtain that all closed Condorcet
domains that contain two completely reversed orders have the structure of a {\em distributive
lattice}. Say that two orders $\overline{R},\underline{R}\in\cR (X)$ are
{\dfn completely reversed} if $\overline{R}\cap\underline{R} =\{ (x,x) \mid x\in X\}$, i.e.,~if
$\overline{R}$ and $\underline{R}$ agree on the ranking of no
pair of distinct alternatives.
\begin{corollary}
\label{clCo-normal-DL} 
Let $\cd$ be a closed Condorcet domain. If $\cd$ contains at least
two completely reversed orders, then it is a distributive
lattice.
\end{corollary}
\noindent{\bf Proof.} By Theorem~\ref{clCo-is-med}, the graph $\Gamma_{\cd}$ associated with $\cd$ is a median graph.
Let $\overline{R},\underline{R}\in\cd$ be two completely
reversed orders, and define the join and meet operations on ${\cd}$ by
$R\vee R' := R^\text{med}(R,\overline{R},R')$ and
$R\wedge R' := R^\text{med}(R,\underline{R},R')$, respectively.  It is easy
(though tedious) to verify that with these operations the space
${\cd}$ has the structure of a distributive
lattice. Another way to prove the corollary is by invoking a well-known result by
\cite{Avann1961} which states that a median graph
is (the covering graph of) a distributive lattice if and only if
it contains two vertices, say $1$ and $0$,
such that every vertex is on a shortest path connecting $1$ and $0$. If
one takes $1$ to be $\overline{R}$ and $0$ to be $\underline{R}$, it is easily seen
using Lemma~\ref{equivalence}(i) that $\Gamma_{\cd}$ satisfies all premises of
Avann's result. % \vspace{-2mm} \\ 
\hspace*{1mm}
\hfill  $\Box$

\medskip

\noindent {\bf Remark.} \cite{Cha-Nem1989} and \cite{Abello91} have shown that all
maximal and connected Condorcet domains that contain at least two completely reversed orders
have the structure of a distributive lattice. This result was generalized by
\cite{DanilovK13} who showed that the connectedness is in fact not needed for the conclusion.
Corollary \ref{clCo-normal-DL} further generalizes this by showing that the condition of maximality
can be substantially weakened to the condition of closedness. In fact, even without the closedness
condition, one would still obtain that every Condorcet domain that contains two completely reversed
orders can be {\em embedded} in a distributive lattice. On the other hand, the condition that the domain
contain at least one pair of completely reversed orders cannot be dropped (an example illustrating
this is given in Figure~\ref{Fig8} below).

%%%%%%%%%%%%%%%%%%%%%%%%%%%%%%%%%%%%%%%%%%%%%%%
%%% SUBSECTION 3.3
%%%%%%%%%%%%%%%%%%%%%%%%%%%%%%%%%%%%%%%%%%%%%%%

\subsection{Every Median Graph is Associated with a Closed Condorcet Domain}

Is every median graph induced by some closed Condorcet domain? The following
result gives an affirmative answer. Interestingly, we will see later in
Section~\ref{maxCond}, that the answer {\color{black} becomes negative}
if we insist on maximality of the
Condorcet domain, i.e., there exist median graphs that cannot be associated with any {\em maximal}
Condorcet domain. 

%%%

\begin{theorem}  
\label{clCo-existence}
For every (finite) median graph $\Gamma = (V,E)$ there exists a closed
Condorcet domain $\cd\subseteq\cR (Y)$ on a finite set of alternatives $Y$ with
$|Y| \leq |V|$ such that $\Gamma_{\cd}$ is isomorphic to $\Gamma$.
\end{theorem}

\noindent
{\bf Proof.} We apply Mulder's theorem.
Since the statement is true for the trivial graph consisting of a single vertex, arguing by
induction, we assume that  the statement is true for all median graphs with  $k$ vertices
or less.  Let $\Gamma'= (V',E')$ be a median graph with $ |V'| = k+1$. By Mulder's
theorem $\Gamma'$ is a convex expansion of some median graph $\Gamma = (V,E)$
relative to convex subsets $W_1$ and $W_2$, where $|V| =\ell\le k$. By induction
there exists a domain $\cd\subseteq \cR(X)$ with $ |X| \le k$ such that
$\Gamma_{\cd}$ is isomorphic to $\Gamma$ with the mapping $R\colon v\mapsto R_v$ 
associating a linear order $R_v\in\cd$ to a vertex $v\in V$.

To obtain a new domain $\cd'$ such that $\Gamma_{\cd'}$ is isomorphic to $\Gamma'$
we clone an arbitrary alternative $x\in X$ and introduce a clone\footnote{We say that $x$ and $y$
are {\em clones} if they are neighbors
in any linear order in the domain, cf.~\cite{Elkind_cloning1,Elkind_cloning2}.} $y\notin X$ of $x$ and denote $X'=X\cup \{y\}$.  The mapping
$R'\colon v\mapsto R'_v$ that associates any vertex $v\in \Gamma'$ to a linear order $R'_v$ over $X'$ will be
constructed as follows. If $v$ is a vertex of $W_1\setminus W_2$, to obtain $R'_v$ we
replace $x$ with $xy$ in $R_v$, placing $x$ higher than $y$, and to obtain $R'_u$ for
$u\in W_2\setminus W_1$ we replace $x$ by $yx$ in $R_u$, placing $y$ higher than $x$. 
Let $v$ now be in $W_1\cap W_2$. In the convex expansion this vertex is split into
$v^1$ and $v^2$. To obtain $R'_ {v^1}$ we clone the linear order $R_v$ replacing
$x$ by $xy$ and to obtain $R'_ {v^2}$ we clone the same linear order $R_v$ replacing
$x$ by $yx$. The number of alternatives has increased by one only, so it is not greater
than $ |V'|=\ell+1\le k+1$.
 
To prove that $\Gamma_{\cd'}$ is isomorphic to $\Gamma'$ we must prove that $R'_u$ and $R'_v$ are neighbors if and only if $u,v\in V'$ are.

First, we need to show that there is no edge between $R'_u$ and $R'_v$ if $v\in W_1\setminus W_2$ and
$u\in W_2\setminus W_1$ since $u$ and $v$ are not neighbors in $\Gamma'$. This follows from the fact that $u$ and $v$ were not neighbors in $\Gamma'$ and hence $R_u$ and $R_v$ were not neighbors in $\Gamma_{\cd}$. Hence there was a linear order
$R_w\in [R_u,R_v]$ between them. In $\cd'$ this linear order will be cloned to
$R'_w$ and, no matter how we place $x$ and $y$ there, we obtain
$R'_w\in [R'_u,R'_v]$ since $R'_u$ and $R'_v$ disagree on $x$ and $y$. Hence $R'_u$ and $R'_v$ are not neighbors in
$\Gamma_{\cd'}$ as well.  Secondly, we have to check that
$R'_ {v^1}$ and $R'_ {v^2}$ are linked by an edge since $v^1$ and $v^2$ are neighbors in $\Gamma'$. This holds because these
orders differ in the ranking of just one pair of alternatives, namely
$x$ and $y$, hence they are neighbors in $\Gamma_{\cd'}$.  If $w\in W_1\cap W_2$ was a neighbor of $v\in W_1\setminus W_2$ in $\Gamma$, then $w^1$ is a neighbor of $v$ in $\Gamma'$ but $w^2$ is not. We then have that $R_w$ is a neighbor of $R_v$ in $\Gamma_{\cd}$. Then $R'_{w^1}$ will obviously be a neighbor of $R'_v$ in $\Gamma_{\cd'}$ while $R'_{w^2}$ will not be a neighbor of $R'_v$ in $\Gamma_{\cd'}$ since $R'_{w^1}$ will be between them. The remaining cases are considered similarly.

To complete the induction step we need to prove that  $\cd'$ is a Condorcet domain.
We use Theorem~\ref{prel_facts}d) for that. We need to consider only triples of elements of
$X'$ that contain $x$ and $y$. Let $\{a,x,y\}$ be such a triple. Since in orders of $\cd'$
elements $x$ and $y$ are always standing together, $a$ can never be in the middle.
Hence conditions of Theorem~\ref{prel_facts}d) are satisfied and $\cd'$ is a Condorcet domain.
\hspace*{1mm}
\hfill \vspace{2ex} $\Box$

\cite{CPS2015} showed that for a median graph with $k$ vertices we might need exactly $k$
alternatives to construct a closed Condorcet domain  that has the given graph
as associated graph, with the star-graph representing the worst-case scenario. Thus,
Theorem~\ref{clCo-existence} cannot be improved  in this respect.
\bigskip

\noindent
{\bf Relation  to the literature.} We conclude this section by commenting on the relation
to \cite{Demange2012}. She introduces profiles of preferences that are
parametrized by a median graph and derives the acyclicity of the majority relation from the
requirement of {\em intermediateness} (cf.~\cite{Grandmont1978}) of preferences. Specifically,
she considers
a median graph with set of vertices $V$ (identified with the set of voters), and a
collection $\{ R_v\in\cR (X) : v\in V\}$ of linear orders and assumes that, for all $v,v',u\in V$,
the following condition is satisfied:
\begin{equation}
\label{Demange-condition}
R_u\in [R_v,R_{v'}] \mbox{ whenever } u \mbox{ is on a shortest path between } v \mbox{ and } v' . 
\end{equation}
Lemma~\ref{equivalence}(i) above shows that, for every median domain $\cd$, the
intermediateness condition \eqref{Demange-condition} is
always satisfied with respect to the geodesic betweenness of the
underlying graph $\Gamma_{\cd}$. Thus, Theorem~\ref{clCo-is-med} above implies
that in fact {\em all} closed Condorcet domains are of the kind
considered by Demange, and demonstrates the existence of a `canonical parametrization' of
these domains (namely by the induced graph $\Gamma_{\cd}$). Finally, \cite{Demange2012}
takes both the (median) graph and the
corresponding profile as given. We show that this does not impose any restriction on the
graph since, by Theorem \ref{clCo-existence}, for every median graph there {\em exists} a
domain satisfying the intermediateness requirement
\eqref{Demange-condition}.\footnote{ In contrast to the present analysis,
\cite{Demange2012} also considers the generalization to weak orders.}

%%%%%%%%%%%%%%%%%%%%%%%%%%%%%%%%%%%%%%%%%%%%%%%
%%%%%%%%%%%%%%%%%%%%%%%%%%%%%%%%%%%%%%%%%%%%%%%
%%%%% SECTION 4
%%%%%%%%%%%%%%%%%%%%%%%%%%%%%%%%%%%%%%%%%%%%%%%
%%%%%%%%%%%%%%%%%%%%%%%%%%%%%%%%%%%%%%%%%%%%%%%

\section{Characterizing and Generalizing the \\ Single-Crossing Property}
\setcounter{equation}{0}
\subsection{The Representative Voter Property}

A {\dfn representative voter} for a given profile of linear orders is a voter, present in this profile,
whose preference coincides with the majority relation of the 
profile. We say that a domain $\cd$ has the {\dfn representative voter property} if in any
odd profile composed of linear orders from~$\cd$  admits a representative voter.

In this subsection, we prove that the single-crossing domains (cf.~Example~1 above) are exactly
the median domains whose associated  graphs are chains. 
From this, we obtain that --- with the exception of the Condorcet domains associated
with the $4$-cycle graph --- the single-crossing domains are the domains that
have the representative voter property.

%%%

\begin{proposition} 
\label{sc-chain}
A domain $\cd\subseteq\cR (X)$ has the single-crossing property if and only if its
associated graph $\Gamma_{\cd}$ is a chain.
\end{proposition}

\noindent {\bf Proof.} It is easily seen that, if $\cd = \{ R_1,\ldots ,R_m \}$ is single-crossing
with respect to the linear order $R_1 > R_2 > \ldots > R_m$, then the interval
$[R_i,R_j]$ for $i<j$ consists of linear orders $R_i,R_{i+1},\ldots, R_j$. Hence, the graph
$\Gamma_{\cd}$ is given by the chain connecting $R_j$ and $R_{j+1}$
by an edge for all $j=1,\ldots,m-1$. 

Conversely, suppose that $\cd=\{ R_1,\ldots,R_m\}$ is such that the associated
graph $\Gamma_{\cd}$ is a chain, say the one that connects $R_j$ with
$R_{j+1}$ by an edge for all $j=1,\ldots,m-1$. We show that $\cd$ is single-crossing
with respect to the linear order $R_1 > R_2 > \ldots > R_m$. By Lemma
\ref{equivalence}(iii), the geodesic betweenness in $\Gamma_{\cd}$ agrees with the Kemeny
betweenness of $\cd$ as a subset of $\cR(X)$, hence, if $R_h,R_l\in\cd$ with $l>h$
are such that $xR_h y$ and $xR_l y$, then we also have $xR_j y$ for all
$j\in\{ h,\ldots, l\}$ and all $x,y\in X$. This implies that $\cd$ has the single-crossing
property relative to the specified order.% \vspace{-2mm} \\ 
\hspace*{1mm}
\hfill \vspace{2ex} $\Box$

%%%

\cite{Roth91} proved the sufficiency of the single-crossing property for the validity
of the representative voter property. The following result shows that the single-crossing property
is in fact `almost' necessary.
\begin{theorem} Let $\cd\subseteq\cR (X)$ be a domain of linear orders on $X$.
Then, $\cd$ has the representative voter property if and only if $\cd$ is either a single-crossing
domain, or $\cd$ is a closed Condorcet domain  with exactly four
elements such that the associated graph $\Gamma_{\cd}$ is a $4$-cycle.
\end{theorem}
\noindent {\bf Proof.} Suppose that $\cd$ has the representative voter property. Evidently,
in this case $\cd$ is a closed Condorcet domain. By Theorem~\ref{clCo-is-med},
the associated graph $\Gamma_{\cd}$ is median, and by Lemma~\ref{equivalence}(i)
the geodesic betweenness in $\Gamma_{\cd}$ agrees with the Kemeny betweenness on $\cd$. We will
show that all vertices of $\Gamma_{\cd}$
have degree at most~$2$.  
Suppose by way of contradiction, that $\Gamma_{\cd}$ contains
a vertex, say $R$, of degree at least~$3$. Consider any profile $\rho = (R_1,R_2,R_3)$
consisting of three distinct neighbors of $R$. The majority relation
corresponding to $\rho$ is the median $R^{\rm med}(R_1, R_2, R_3)$. Since
the median graph $\Gamma_{\cd}$ does not have 3-cycles, $R_i$ and $R_j$ cannot
be neighbors for all
distinct $i,j\in\{ 1,2,3\}$, hence $R^{\rm med}(R_1, R_2, R_3) = R$.
Since $R$ is not an element of $\{R_1,R_2,R_3\}$ the representative voter
property is violated, a contradiction.
Since $\Gamma_{\cd}$ is always connected (as a graph),
the absence of vertices of degree 3 or more implies that $\Gamma_{\cd}$ is either a cycle or
a chain. It is well-known that
among all cycles, only $4$-cycles are a median graphs.\footnote{\label{cycle-fn}Clearly, a 
$3$-cycle is not a median graph. Moreover, for any $k\geq 5$, one can find three vertices
$v_1,v_2,v_3$ on a $k$-cycle such that the three shortest paths between any pair
from $v_1,v_2,v_3$ cover the entire cycle; this implies that no vertex can simultaneously
lie on all three shortest paths, i.e.,~that the triple $v_1,v_2,v_3$ does not admit a median.}
On the other hand, if $\Gamma_{\cd}$ is a chain, then $\cd$ has the single-crossing property
by Proposition~\ref{sc-chain}.

To prove the converse, suppose first that $\cd$ is a single-crossing domain.
Then $\Gamma_{\cd}$ is a chain by Proposition~\ref{sc-chain} and, evidently,
the preference of the median voter in any odd profile coincides with the corresponding
majority relation; this is Rothstein's theorem (\cite{Roth91}). On the other hand,
consider any odd profile over a domain $\cd$ such that the induced graph is
a $4$-cycle. In that case, the representative voter property holds trivially if the
profile contains all four different orders; if it contains at most three different orders,
the representative voter property follows as in the case of a chain with at most three
elements.% \vspace{-2mm} \\ 
\hspace*{1mm} \hfill $\Box$

%%%

\subsection{Generalizing the Single-Crossing Property to Trees}

The classical single-crossing property  \citep{Mirr71, GansSmart1996} naturally generalizes to trees.
We rephrase the definition suggested by~\cite{Kung2014} and~\cite{CPS2015} as follows.

\begin{definition}
A domain $\cd\subseteq\cR (X)$ is {\dfn single-crossing  with respect to the}  tree
$T=(V,E)$ if $|\cd|=|V|$ and the linear orders of $\cd$ can be parametrized by the vertices of $T$ so
that the intermediateness condition \eqref{Demange-condition} is satisfied.  Moreover,
we say that a domain is {\dfn generalized single-crossing} if it is single-crossing with respect to
some tree.
\end{definition}

%Every such domain is a Condorcet domain \citep{Demange2012,CPS2015}.  \par\medskip

As is easily verified, a domain $\cd$ is generalized single-crossing if and only if its associated
graph $\Gamma_{\cd}$ is a tree. Indeed, if a domain $\cd$ is single-crossing with respect
to the tree $T$, then $\Gamma_{\cd}=T$. Conversely, if $\Gamma_{\cd}$ is a tree, then 
by Lemma~\ref{equivalence}(iii), the preferences of $\cd$ are intermediate on $\Gamma_{\cd}$
so $\cd$ is generalized single-crossing. In particular, every generalized single-crossing
domain is a closed Condorcet domain by Lemma \ref{equivalence}(iii)
(cf.~\cite{Demange2012,CPS2015}).

 To justify our terminology of `generalized {\em single-}crossingness,' let
$\cd\subseteq \cR(X)$
be a generalized single-crossing domain and suppose not all orders
in $\cd$ agree on the ranking of $x,y\in X$, then the sets $\cvd_{xy} := \{ R\in\cd \mid xRy\}$
and $\cvd_{yx} := \{ R\in\cd \mid yRx\}$ are both non-empty. In the tree $\Gamma_\cd$, there
will be a {\em unique} edge $QR$ such that
$Q\in\cvd_{xy}$ and $R\in\cvd_{yx}$; moreover, both sets $\cvd_{xy}$ and  $\cvd_{yx}$ are
convex, and $\cvd_{xy}\cup \cvd_{yx}=\cd$. This justifies the terminology 

\par\medskip

\noindent {\bf Example 3 (The single-crossing property on a tree).}
Consider domain $\cd$ on the set $\{ a,b,c,d\}$ consisting of 
four orders: $\hat{R} = abcd$, $R_1 = acbd$, $R_2 = abdc$, and $R_3 = bacd$.
As is easily seen, $\cd$ is a closed Condorcet domain. The  associated graph
$\Gamma_{\cd}$ connects $\hat{R}$ with each of the
other three orders by an edge and the graph has no other edges (cf.~Figure~4). Thus,
$\Gamma_{\cd}$ is a tree and $\hat{R}$ is the median order of any
triple of distinct elements of $\cd$.  

\begin{figure}[h]
\label{Fig4}
\centering
\includegraphics[width=8cm]{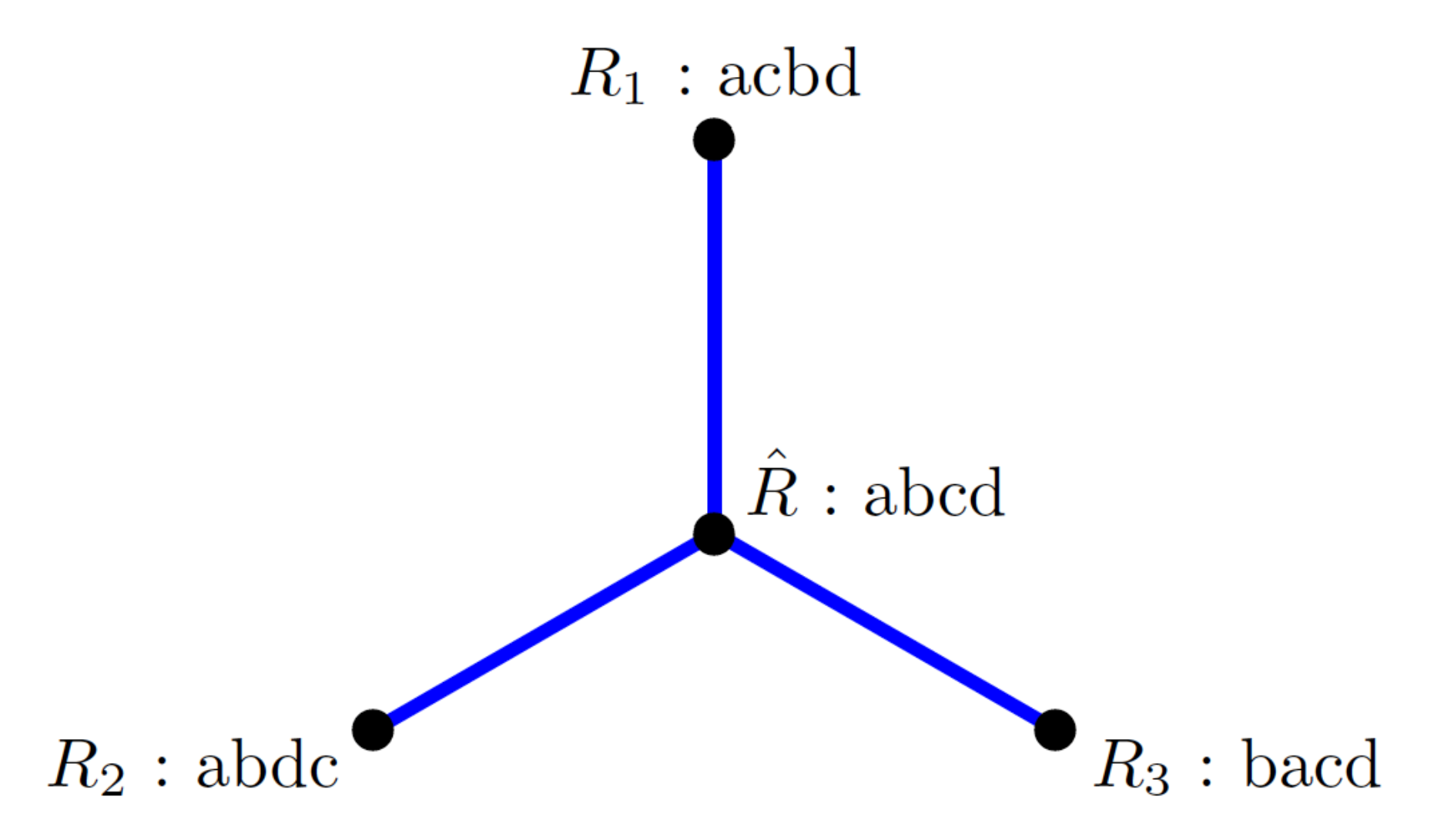}
\caption{A generalized single-crossing domain}
\end{figure}

%
%\vbox{
%\begin{center} 
%\hspace*{0mm}\includegraphics[width=8cm]{Fig-SingleCrossTree}
%\end{center}
%
%\begin{center}
%{\em Figure 4: The single-crossing property on a tree}
%\vspace{3mm}
%\end{center}
%}
%FI

The following result characterizes the single-crossing property and its generalization to
trees directly in terms of the structure of the underlying domain, i.e., without explicit reference
to the associated graph.

%%%

\begin{proposition} {\bf a)} A domain $\cd$ has the single-crossing property with respect
to some linear order on $\cd$ if and only if, for all $x,y,z,w\in X$ such that each of the sets
$\cvd_{xy}$, $\cvd_{yx}$, $\cvd_{zw}$, and $\cvd_{wz}$ is non-empty, we have:
\begin{equation} \label{eq-single-cross}
\cvd_{xy}\subseteq \cvd_{zw}\mbox{ or }\cvd_{xy}\subseteq \cvd_{wz}.
\end{equation}
{\bf b)} Let ${\mathbb F}=\{\cvdxy \mid x,y\in X\}$. A domain $\cd$ has the generalized
single-crossing property if and only if ${\mathbb F}$ has the
Helly property and,  for all distinct $x,y,z,w\in X$, at least one of the following
four sets is empty:
\begin{equation} \label{eq-gen-single-cross}
\cvd_{xy}\cap\cvd_{zw},\ \cvd_{xy}\cap\cvd_{wz},\
\cvd_{yx}\cap\cvd_{zw},\ \cvd_{yx}\cap\cvd_{wz}.
\end{equation}
\end{proposition}

\noindent
{\bf Proof. a)} Evidently, every single-crossing domain  satisfies \eqref{eq-single-cross}.
To prove the converse, we notice that this condition allows one to order the family
of all sets of the form $\cvd_{xy}$ so that, for an appropriate sequence of pairs
of alternatives $(a_1,b_1), \ldots , (a_m,b_m)$,
\[
\cvd_{a_1 b_1}\subseteq \cvd_{a_2 b_2}\subseteq \ldots \subseteq \cvd_{a_m b_m}
\ \mbox{and} \
\cvd_{b_1 a_1}\supseteq \cvd_{b_2 a_2}\supseteq \ldots \supseteq \cvd_{b_m a_m}.
\]
It is now easily verified that $\cd$ has the single-crossing property with respect to any linear
order of the members of $\cd$ which lists the elements of $\cvd_{a_1 b_1}$ first, then the
elements of $\cvd_{a_2 b_2}\setminus \cvd_{a_{1} b_{1}}$ second, and further lists elements
of $\cvd_{a_j b_j}\setminus \cvd_{a_{j-1} b_{j-1}}$ after listing $\cvd_{a_{j-1} b_{j-1}}$.
\medskip

{\bf b)} Suppose that $\cd$ has the generalized single-crossing property. Then,
$\Gamma_{\cd}$ is a tree. By Lemma \ref{equivalence}(iii), the betweenness
on $\cd$ coincides with the induced geodesic betweenness on $\Gamma_{\cd}$.
In particular, $\cd$ is a median domain, and as such satisfies the Helly property for
convex sets, and all elements of ${\mathbb F}$ are convex. To
verify \eqref{eq-gen-single-cross}, assume by contradiction that
\[
P\in \cvd_{xy}\cap\cvd_{zw},\ Q\in \cvd_{xy}\cap\cvd_{wz},\
R\in \cvd_{yx}\cap\cvd_{zw},\ S\in \cvd_{yx}\cap\cvd_{wz},
\]
for $P,Q,R,S\in\cd$. Consider a shortest path between $P$ and $R$ and a shortest path
between $Q$ and $S$. The first path lies entirely in $\cvd_{zw}$ and the second one lies
entirely in $\cvd_{wz}$; in particular, they do not intersect. But on each of them there is
a switch from $xy$ to $yx$, and thus there exist two pairs of neighboring orders such that one
of them is in $\cvd_{xy}$ and the other one in $\cvd_{yx}$. This contradicts the generalized
single-crossing property.

Conversely, suppose that a domain $\cd$ satisfies the Helly property and condition
(\ref{eq-gen-single-cross}). Proposition~\ref{helly-prop} then implies that $\cd$ is a median domain,
hence by Lemma \ref{equivalence}(i), the betweenness in $\cd$ coincides with the geodesic
betweenness in $\Gamma_{\cd}$. By (\ref{eq-gen-single-cross}), $\Gamma_{\cd}$
is acyclic, hence a tree, which implies the generalized single-crossing property, as desired.
%\vspace{-1mm} \\  
\hspace*{1mm} \hfill $\Box$

%%%%%%%%%%%%%%%%%%%%%%%%%%%%%%%%%%%%%%%%%%%%%%%
%%%%%%%%%%%%%%%%%%%%%%%%%%%%%%%%%%%%%%%%%%%%%%%
%%% SECTION 5
%%%%%%%%%%%%%%%%%%%%%%%%%%%%%%%%%%%%%%%%%%%%%%%
%%%%%%%%%%%%%%%%%%%%%%%%%%%%%%%%%%%%%%%%%%%%%%%

%
\section{Maximal Condorcet Domains Revisited \label{maxCond}}
\setcounter{equation}{0}
The problem of characterizing {\em maximal} Condorcet domains has received
considerable attention in the literature. Since every maximal Condorcet domain is
closed, our analysis allows one to address this problem by studying the structure
of the associated median graph. In particular, we will show that, if the median graph
associated with a maximal domain is acyclic, then it must, in fact, be a chain, i.e.,
among all trees only chains can be associated with maximal Condorcet domains.
On the other hand, it is well-known that single-crossing domains are in general
not maximal as Condorcet domains (see, e.g., \cite{mon:survey}).
However, occasionally this happens, and we provide a
simple necessary and sufficient
condition for a single-crossing domain to be a maximal Condorcet domain. Finally, we demonstrate
by means of an example that even though the induced median graph of
a maximal Condorcet domain is never a tree different from a chain, it does not need
to have the structure of a distributive lattice either (in view of
Corollary~\ref{clCo-normal-DL} above, such
maximal Condorcet domain cannot
contain two completely reversed orders).

\begin{figure}[h]
\label{Fig5}
\centering
\includegraphics[width=10cm]{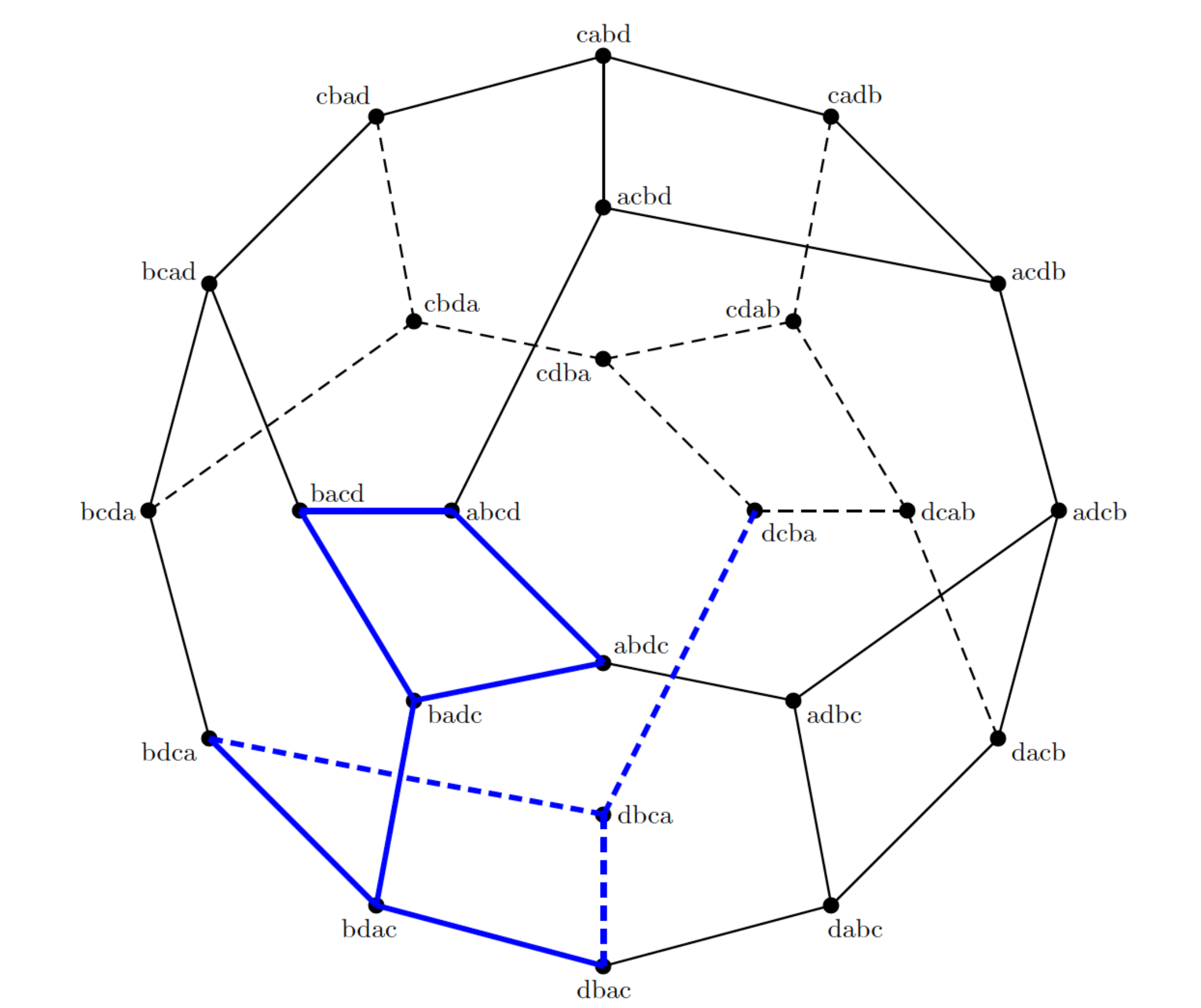}
\caption{A connected maximum Condorcet domain on $\{ a,b,c,d\}$}
\end{figure}

If $X$ has three elements, all maximal Condorcet domains have four elements and are either of the
type shown in Figure~\ref{Fig2} (connected and single-crossing) or of the type shown to the left of
Figure~\ref{Fig3} (not connected and cyclic). If $ |X| > 3$, not all maximal Condorcet domains on $X$
%$\cd\subseteq\cR (X)$ 
have the same cardinality. Call a domain a
{\dfn maximum} Condorcet domain on~$X$
if it achieves the largest cardinality among all Condorcet domains on~$X$. Evidently, every
maximum Condorcet domain is also maximal (and therefore closed), but not vice versa.
If $X$ has four elements, the cardinality
of a maximum Condorcet domain is known to be~$9$, see Figure~\ref{Fig5} for a domain
attaining this number.

Note that the domain shown in Figure~\ref{Fig5} is connected. Indeed, most of the known general
results about the size and the structure of maximal Condorcet domains
pertain to {\em connected} Condorcet domains; for instance, the domains described in
\cite{Abello91,Cha-Nem1989,DKK:2012} and \cite{GR:2008} are all connected.
However, not all maximal Condorcet domains are connected, as we have already seen
above. %(cf.~Figure~\ref{Fig3} to the left and recall that even though the induced graph $\Gamma_{\cd}$
%is always a connected graph, the underyling domain need not be connected as a domain).
Figure 6 shows a maximal Condorcet domain
with $8$ elements on $\{ a,b,c,d\}$ that is not connected; in fact, it belongs to the
class of the so-called `symmetric' Condorcet domains studied in \cite{DanilovK13};
a domain is called {\dfn symmetric} if it contains with any order also its completely reversed
order.

%\vspace{3ex}
%\vspace*{-25mm}
%\begin{center} 
%\hspace*{0mm}\includegraphics[width=10cm]{Fig-MaximumCond}
%\vspace*{2ex}
%\end{center}
%
%\begin{center}
%{\em Figure 5: A connected maximum Condorcet domain on $\{ a,b,c,d\}$}
%\vspace{2ex}
%\end{center}
%

\begin{figure}[h]
\label{Fig6}
\centering
\includegraphics[width=10cm]{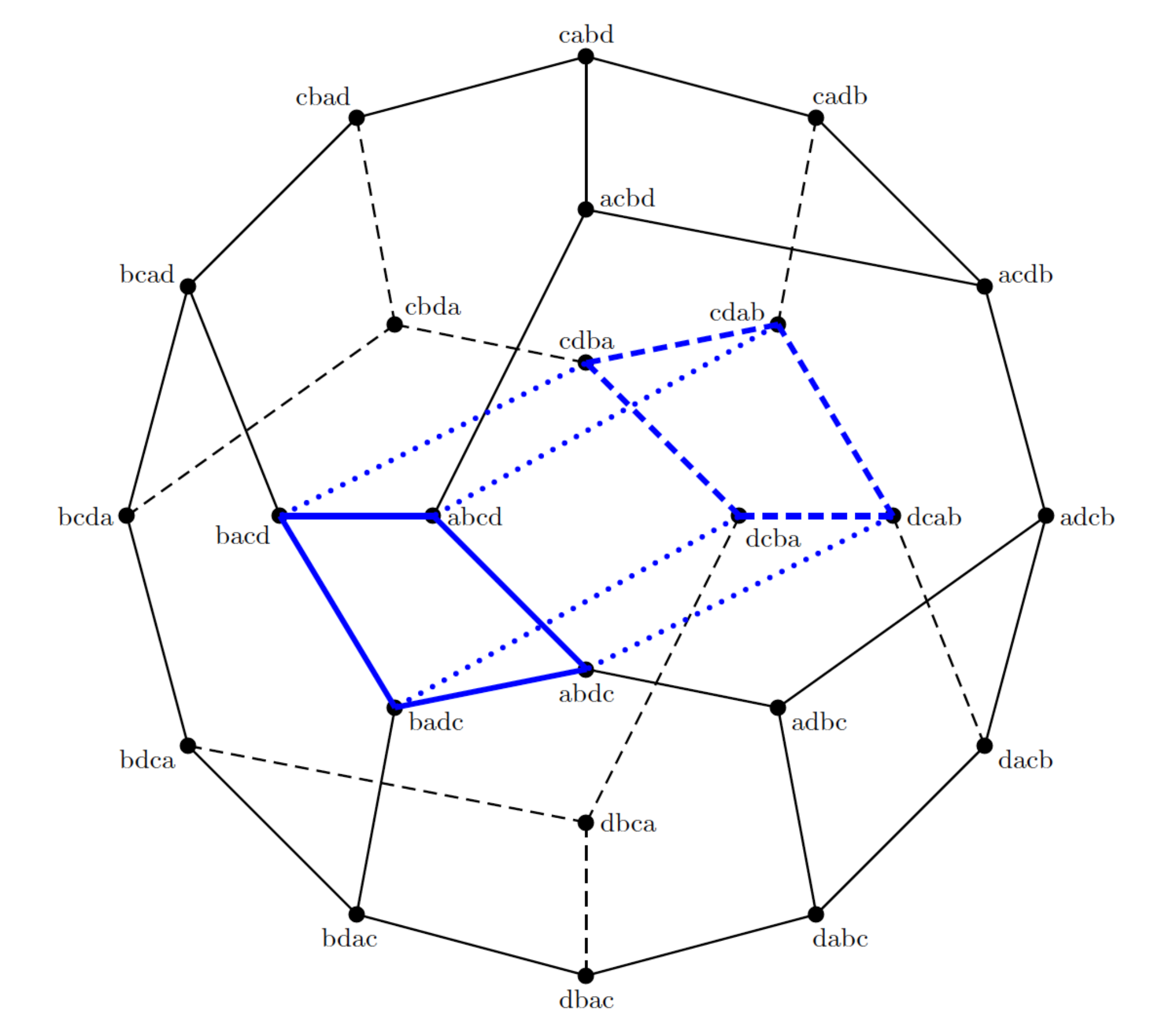}
\caption{A non-connected maximal Condorcet domain on $\{ a,b,c,d\}$}
\end{figure}
%
%\vspace{2ex} 
%\begin{center} 
%\hspace*{0mm}\includegraphics[width=10cm]{Fig-MaxCube}
%\vspace{2ex} 
%\end{center}
%
%\begin{center}
%{\em Figure 6: A non-connected maximal Condorcet domain on $\{ a,b,c,d\}$}
%\vspace{2ex} 
%\end{center}
%

All examples of maximal Condorcet domains that we have given so far are either
single-crossing domains, or their associated graph contains
at least one $4$-cycle.\footnote{It is well-known that every median graph that is not
a tree contains at least one $4$-cycle; this can be verified as in Footnote~\ref{cycle-fn}
above. Note that median graphs can contain also larger cycles (cf.Figure~\ref{Fig6})
but these can never be
`minimal' (in an appropriate sense); every minimal cycle of any median graph
is a $4$-cycle.} As we will see, this is a general feature of all
maximal Condorcet domains. To show this we need to do some preliminary work.

 \begin{lemma}
\label{not_maximal}
Suppose $\cd$ is generalized single-crossing and  $QR$ is an edge in $\Gamma_\cd$. Suppose $P\in [Q,R]$ is different from $Q$ and $R$.  Then $\cd'=\cd\cup \{P\}$ is also a generalised single-crossing domain and, in particular, $\cd$ is not a maximal Condorcet domain.
\end{lemma} 
{\bf Proof.} We add $P$ to the graph splitting the edge $QR$, placing $P$ in the middle.
The new graph will be still a tree and the preferences of $\cd'=\cd\cup \{P\}$ will be intermediate
with respect to this tree, hence generalized single-crossing. In particular, $\cd'=\cd\cup \{P\}$ is
a median domain and hence a closed Condorcet domain by Theorem~\ref{cc-iff-med}. \hfill $\Box$

\par\medskip
\noindent
We note that the lemma is specific to trees; if $\Gamma_\cd$ is a rectangle, the corresponding
statement would not hold as the right-hand side of Figure~\ref{Fig3} shows. 

\begin{corollary}
\label{gsc-connected}
Any maximal generalized single-crossing domain on $X$ is connected, i.e., is a subgraph of the permutohedron.
\end{corollary}
%%%

Now we can prove the first of the two main results of this section.

\begin{theorem}
\label{maximal-tree-is-chain}
Let $\cd$ be a maximal Condorcet domain. If $\Gamma_{\cd}$ is a tree, it is, in fact, a chain.
\end{theorem}

\noindent
{\bf Proof.} Let $\cd$ be a maximal Condorcet domain, and assume that $\Gamma_{\cd}$
is a tree but not a chain. Then there exists a vertex $R$ in $\Gamma_{\cd}$ of degree at least~$3$.   
Consider now any three neighbors of $R$ in $\Gamma_{\cd}$, say $R_1$, $R_2$ and $R_3$. 
Since, by Corollary~\ref{gsc-connected}, $\cd$ is connected, there are
 three distinct ordered pairs $(x_i,y_i)$, $i=1,2,3$, of  alternatives such that 
 $R_i=R\setminus \{(x_i,y_i)\}\cup \{(y_i,x_i)\}$. We will say that $R_i$ is obtained from $R$ by {\dfn switching} the pair of adjacent alternatives $(x_i,y_i)$.
Moreover, since in every pair $(x_i,y_i)$, $i=1,2,3$, the alternatives are adjacent in $R$, there must exist at least two pairs that have no alternative
in common, say $\{ x_1,y_1\} \cap \{ x_2,y_2\} = \emptyset$. 

Now let $R'$ be the order that coincides with $R$ except that both pairs $(x_1,y_1)$ and
$(x_2,y_2)$ in $R'$ are switched, i.e., $y_1 R' x_1$ and
$y_2 R' x_2$, and consider the domain $\cd\cup\{ R'\}$.
Since $x_1,y_1$ and $x_2,y_2$ are neighbors in each of the orders $R$, $R_1$, $R_2$, $R'$,
for every three alternatives $\{a,b,c\}$ no new order among them appears in $R'$ which has not
yet occurred in $R$, $R_1$, or $R_2$. Hence, by Theorem~\ref{prel_facts}d), $\cd\cup\{ R'\}$
is a Condorcet domain. By the maximality of $\cd$, this implies $R'\in\cd$. But in this case,
the graph $\Gamma_{\cd}$ evidently contains the $4$-cycle $\{ R,R_1,R',R_2\}$,
contradicting the assumed acyclicity $\Gamma_{\cd}$. Hence, there cannot exist a vertex of
degree $3$ or larger, i.e., $\Gamma_{\cd}$ is a chain.
%\vspace{-1mm} \\  
\hspace*{1mm} \hfill \vspace{2ex} $\Box$

 Figure~4 (cf.~Example~3 above) illustrates the proof of Theorem
\ref{maximal-tree-is-chain}: one easily
verifies that the order $badc$ can be added to the depicted domain $\cd$, creating
a $4$-cycle in the associated graph $\Gamma_{\cd}$; in particular, $\cd$
is not maximal.
\par\medskip
%%%

As we will see some -- but by far not all -- maximal single-crossing domains are also maximal as Condorcet domains. Figure~\ref{Fig5}, however, shows instances of non-maximality:  the depicted maximal Condorcet domain contains four maximal single-crossing domains as proper subdomains.
The next result characterizes exactly when a maximal single-crossing domain is  also a maximal Condorcet domain. \par\medskip

To formulate the result, we need the following definitions. Let $\cd\subseteq\cR(X)$
be a maximal single-crossing domain and $|X|=n$. Then we know that the associated graph
$\Gamma_\cd$ is linear. By Corollary~\ref{gsc-connected}, it is a subgraph of the
permutohedron. Let us enumerate 
orders of $\cd$ so that the edges of $\Gamma_\cd$ are $R_1R_2$, $R_2R_3$,
$\ldots,$ $R_{m-1}R_m$. Then the sequence $\{\row Rm\}$ will be called a
{\dfn maximal chain}. Due to Lemma~\ref{not_maximal},
each edge $R_iR_{i+1}$ in $\Gamma_\cd$ has a unique switching pair of alternatives $(x_i,y_i)$
which are adjacent in $R_i$. Due to the (standard) single-crossing condition, once switched the pair is
never switched back, so for a two-element subset of alternatives $\{x,y\}$ we cannot
have both pairs $(x,y)$ and $(y,x)$ as switching pairs. Moreover, $R_1$ and $R_m$
are completely reversed. Indeed, if $R_m$ is not $\underline{R_1}$, which is completely
reversed $R_1$, then we can add $\underline{R_1}$ at the end of the sequence $\{\row Rm\}$
and add another vertex to our graph in contradiction to maximality  of $\cd$. Hence,
$m=\frac{(n-1)n}{2}+1$, so, in particular, all maximal chains are of the same length.
\par\smallskip

Without loss of generality we may assume that $X=\{1,2,\ldots, n\}$ and
$\overline{R}=R_1=1\,2\ldots n-1\, n$ and $\underline{R}=R_m=n\, {n-1}\ldots 2\, 1$. 
Traveling along the maximal chain $\row Rm$ from $\overline{R}$ in direction of $\underline{R}$ alternatives in every 
2-element subset $\{x,y\}$, where $x<y$ switch exactly once from $xy$ to $yx$. 
That is, for each pair $(x,y)\in X\times X$
there exists exactly one edge $(R_j,R_{j+1})$ such that $xR_1 y,\ldots, xR_j y$
and $yR_{j+1}x,\ldots, yR_m x$.
%We will refer to such pair as the
%{\dfn switching pair} corresponding to the edge of $\Gamma_{\cd}$. 
We say that this
maximal chain satisfies the {\dfn pairwise concatenation condition}
if the switching pairs corresponding to any two adjacent edges have one alternative
in common.  Formally, a maximal chain $\{\row Rm\}$ with switching pairs $(x_j,y_{j})$ of edges $R_jR_{j+1}$, $j=1,...,m-1$, satisfies pairwise concatenation condition if
\begin{equation} 
\label{eq:pcc}
\{ x_j,y_j\} \cap \{ x_{j+1},y_{j+1}\} \neq\emptyset 
\end{equation}
for all $j=1,...,m-1$.
 Note
that the intersection in~\eqref{eq:pcc} then  has exactly one element. We note that the
sequence of switching pairs $(x_1,y_1), \ldots, (x_{m-1},y_{m-1})$ determines the maximal
chain $\{\row Rm\}$ uniquely.

To illustrate
the pairwise concatenation condition consider Figure~7 which depicts a single-crossing domain
on $X = \{ a,b,c,d\}$ with $m=7$. The maximal chain for this domain is:
\[
 R_1 = abcd,\ R_2=acbd,\ R_3=acdb,\ R_4=adcb,\ R_5=dacb,\ R_6=dcab,\  R_7 = dcba, 
 \]
 with the sequence of switching pairs $(b,c), (b,d), (c,d), (a,d), (a,c), (a,b)$, which obviously satisfies the pairwise concatenation
condition.

\begin{figure}
\label{Fig7}
\centering
\includegraphics[width=10cm]{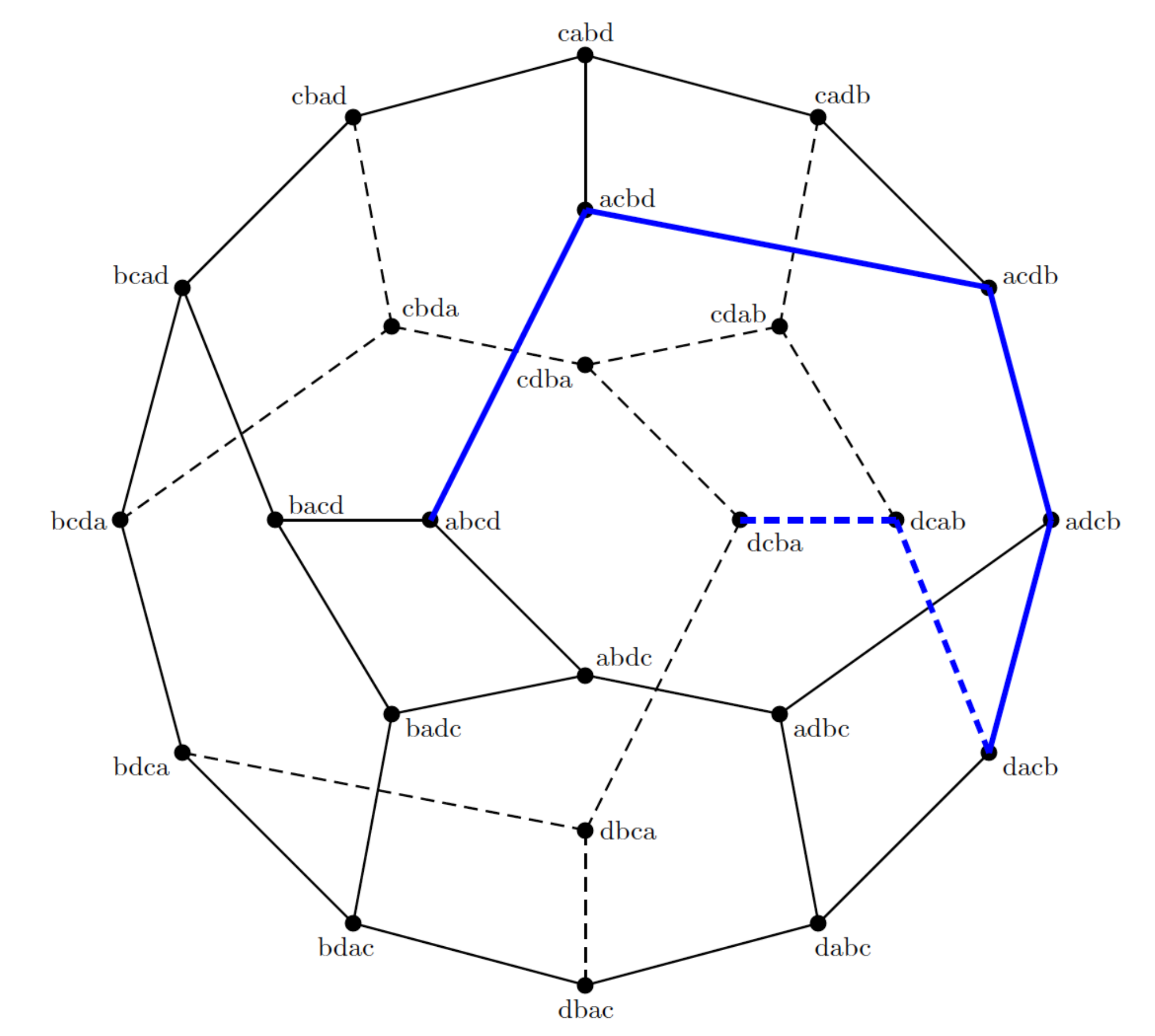}
\caption{A maximal chain constituting a maximal Condorcet domain on $\{ a,b,c,d\}$}
\end{figure}

We will now show that the pairwise concatenation condition is necessary and
sufficient for a maximal single-crossing domain to be a maximal Condorcet domain.\footnote{As we
learned after completion of the first version of the present paper, Bernard Monjardet
found the same condition in his unpublished lecture notes \cite{Monjardet2007}
(without proof).} %But first we need to prove the following

\begin{definition}[\cite{GR:2008}] 
Two maximal chains are {\dfn equivalent} if their respective sequences of switching pairs can be transformed one to another by swapping adjacent non-intersecting pairs.
\end{definition}

 It is easily verified that this notion indeed defines an
equivalence relation on the class of all maximal chains. To illustrate the definition,
consider the maximum Condorcet domain $\cd$ on
the set $X=\{ a,b,c,d\}$ depicted in Figure \ref{Fig5} above. The following sequence
is one of its maximal chains:
\[
\cc_1 := \left< \ abcd,\ abdc,\ badc,\ bdac,\ dbac,\ dbca,\ dcba \ \right>
\]
with the corresponding sequence of switching pairs
\begin{equation} \label{switch_pair}
(d,c), (a,b), (a,d), (b,d), (a,c), (b,c).
\end{equation}
If one swaps the first two pairs in this sequence, i.e.~$(d,c)$ and $(a,c)$, one obtains the
sequence of switching
pairs corresponding to the equivalent maximal chain:
\[
\cc_2 := \left< \ abcd,\ bacd,\ badc,\ bdac,\ dbac,\ dbca,\ dcba \ \right> ,
\]
which is also contained in the maximum domain $\cd$. If in the sequence \eqref{switch_pair},
one switches the fourth and fifth pair, i.e.~$(b,d)$ and $(a,c)$, one obtains another
equivalent maximal chain of $\cd$:
\[
\cc_3 := \left< \ abcd,\ abdc,\ badc,\ bdac,\ bdca,\ dbca,\ dcba \ \right> .
\]
Finally, if one switches both pairs in \eqref{switch_pair}, first $(d,c)$ against $(a,b)$ and then
$(b,d)$ against $(a,c)$, one obtains the equivalent maximal chain:
\[
\cc_4 := \left< \ abcd,\ bacd,\ badc,\ bdac,\ bdca,\ dbca,\ dcba \ \right> .
\]
Evidently, the maximum Condorcet domain $\cd$ is the union of the pairwise equivalent
maximal chains $\cc_1$, $\cc_2$, $\cc_3$, $\cc_4$ that it contains. The following result due
to \cite[Th.~2]{GR:2008} shows that this
is a general property of all maximal Condorcet domains that contain at least one maximal
chain.

\begin{theorem}[\cite{GR:2008}]
\label{GR-theorem}
Let $\cd$ be a single-crossing domain and let $\left< \row Rm\right>$ be its corresponding
maximal chain.
Then the maximal Condorcet domain containing $\cd$ consists of all linear orders of all
maximal chains equivalent to $\left< \row Rm \right>$. 
\end{theorem}

We can now state the second main result of this section.

\begin{theorem}
\label{max-single-cross}
A maximal single-crossing domain $\cd\subseteq\cR (X)$ is a maximal Condorcet
domain if and only if  the  maximal chain corresponding to $\cd$ satisfies the pairwise concatenation condition.  
\end{theorem}

%%%
 
%\begin{theorem}
%\label{max-single-cross}
%A maximal single-crossing domain $\cd\subseteq\cR (X)$ is a maximal Condorcet
%domain if and only if (i) $\cd$ is connected, (ii) $\cd$ contains two completely reversed
%orders, and (iii) the corresponding maximal chain satisfies the pairwise concatenation condition.  
%\end{theorem}

%
%\vspace{-2mm}
%\begin{center} 
%\hspace*{0mm}\includegraphics[width=10cm]{Fig-MaxChain}
%\vspace{1mm} 
%\end{center}

%\begin{center}
%{\em Figure 7: A maximal chain constituting a maximal Condorcet domain on $\{ a,b,c,d\}$}
%\vspace{2ex} 
%\end{center}
%

\noindent
{\bf Proof.} Let $\cd$  be a maximal single-crossing domain and $\{\row Rm\}$ be the
corresponding maximal chain.
To verify the necessity of the pairwise concatenation condition suppose, by contraposition, that for two
consecutive switching pairs $(x_j,y_j)$ and $(x_{j+1},y_{j+1})$ one has
$\{ x_j,y_j\}\cap\{ x_{j+1},y_{j+1}\}=\emptyset$. Without loss of generality assume that
$x_jR_jy_j$ and $y_jR_{j+1}x_j$, as well as $x_{j+1}R_{j+1}y_{j+1}$ and
$y_{j+1}R_{j+2}x_{j+1}$. Since $\{ x_j,y_j\}\cap\{ x_{j+1},y_{j+1}\}=\emptyset$,
both pairs of alternatives $(x_j,y_j)$ and $(x_{j+1},y_{j+1})$ must be adjacent
in all three orderings $R_j,R_{j+1},R_{j+2}$, respectively. In particular, these
three orders agree in the ranking of all pairs except for the two pairs $(x_j,y_j)$ and
$(x_{j+1},y_{j+1})$, respectively. Consider the order $R'$ that has $x_jR'y_j$
and $y_{j+1}R'x_{j+1}$ and agrees with each of $R_j,R_{j+1},R_{j+2}$ in the ranking of
all other pairs. Clearly, $R'\not\in\cd$, but the restriction of $R'$ to any triple of distinct
alternatives coincides with the restriction of $R_j$, $R_{j+1}$, or $R_{j+2}$ to this
triple. Hence, $\cd\cup\{ R'\}$ is a Condorcet domain by Theorem \ref{prel_facts}d),
and thus $\cd$ is not maximal. This demonstrates the necessity of the pairwise
concatenation condition. 

The sufficiency  of the pairwise concatenation condition follows from
Theorem~\ref{GR-theorem}. Indeed, it is easily seen that a maximal chain
satisfying the pairwise concatenation condition has no other maximal chain equivalent to it,
i.e.~it forms an equivalence class of its own. By Theorem~\ref{GR-theorem},
the corresponding maximal single-crossing domain is a maximal Condorcet domain. \hfill $\Box$ 
\bigskip

%Geometrically, the orders $R_j,R_{j+1},R_{j+2},R'$ form a $4$-cycle. For illustration,
%consider the orders $R_j = abcd$, $R_{j+1} = bacd$, $R_{j+2} = badc$ in Fig.~5 above,
%with the corresponding switching pairs $(x_j,y_j)=(a,b)$ and $(x_{j+1},y_{j+1})=(c,d)$.
%The order $R'$ constructed above is given here by $abdc$; it arises by switching
%both pairs $(a,b)$ and $(c,d)$ in $R_{j+1}$ simultaneously, and thus completes the
%$4$-cycle.

\begin{figure}
\centering
\includegraphics[width=9cm]{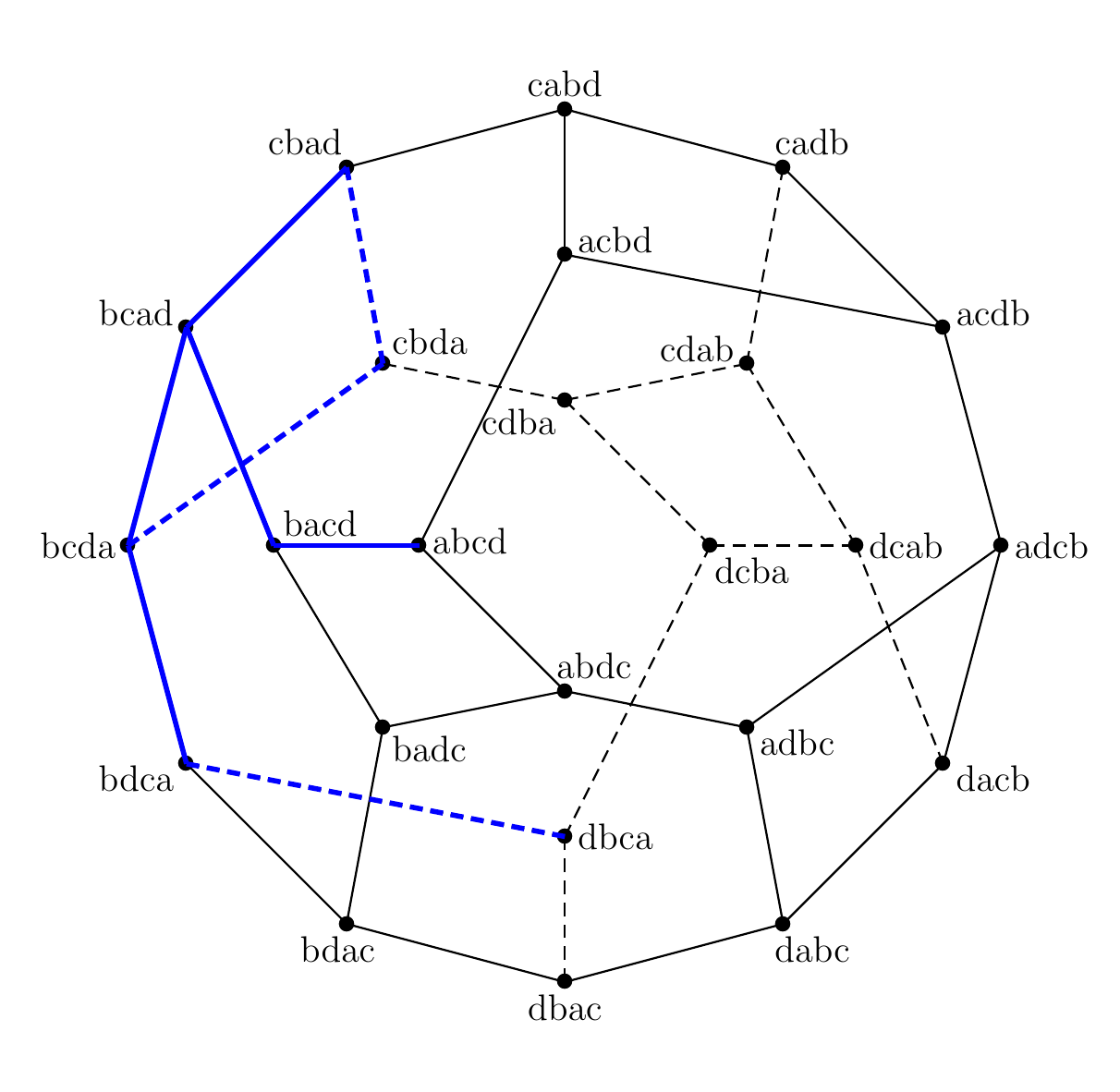}
\caption{A maximal Condorcet domain with no pair of completely
reversed orders}
\label{Fig8}
\end{figure}

In light of Theorem~\ref{maximal-tree-is-chain} and Corollary \ref{clCo-normal-DL}
above, a natural conjecture is that all maximal Condorcet domains are distributive lattices.
However, this turns out not to be true, as the example in the following figure shows.
The domain depicted in Figure~\ref{Fig8} is connected and maximal but not a distributive lattice;
in view of Corollary \ref{clCo-normal-DL} it clearly cannot
contain a pair of completely reversed orders.

%%%

%\vspace*{-2mm}
%\begin{center} 
%\includegraphics[width=9cm]{Figure_8a}
%\vspace{1mm} \\
%{ {\em Figure 8: A maximal Condorcet domain with no pair of completely reversed orders}}
%\vspace{2mm}
%\end{center}

%%%

It seems to be an interesting and worthwhile subject for future research to characterize
the class of median graphs that can arise as the associated graphs of {\em maximal}
Condorcet domains. However, this approach has to be complemented by other
considerations as well,
since the graph alone does in general not contain information about maximality. For
instance, we have already seen that chains may be the associated graphs of
maximal Condorcet domains only if their length is one greater than a triangular
number, i.e., $\frac{1}{2}(n-1)n +1$. 
Another example that illustrates this point quite
drastically is the $4$-cycle. We have already seen above that even non-median
domains can induce a $4$-cycle, albeit the betweenness structures of the domain
and the $4$-cycle will not be the same (see the domain in the
middle of Figure~\ref{Fig3}). But even if we insist
on them being isomorphic, thus ensuring that the corresponding domain is a closed
Condorcet domain, the $4$-cycle may or may not yield a {\em maximal} Condorcet domain.
For instance, all $4$-cycles in Figures 5, 6 and 8 evidently do not correspond to maximal
Condorcet domains since they form proper subdomains. But  Figure 9 shows that the
$4$-cycle can also be the associated graph of a maximal domain in
$\cR (\{ a,b,c,d\} )$.  Such maximal domain of size $4$ can be found in
$\cR(X)$ for any cardinality of $X$, as shown by \cite{DanilovK13}.

%%%

\begin{figure}[h]
\centering
\includegraphics[width=9cm]{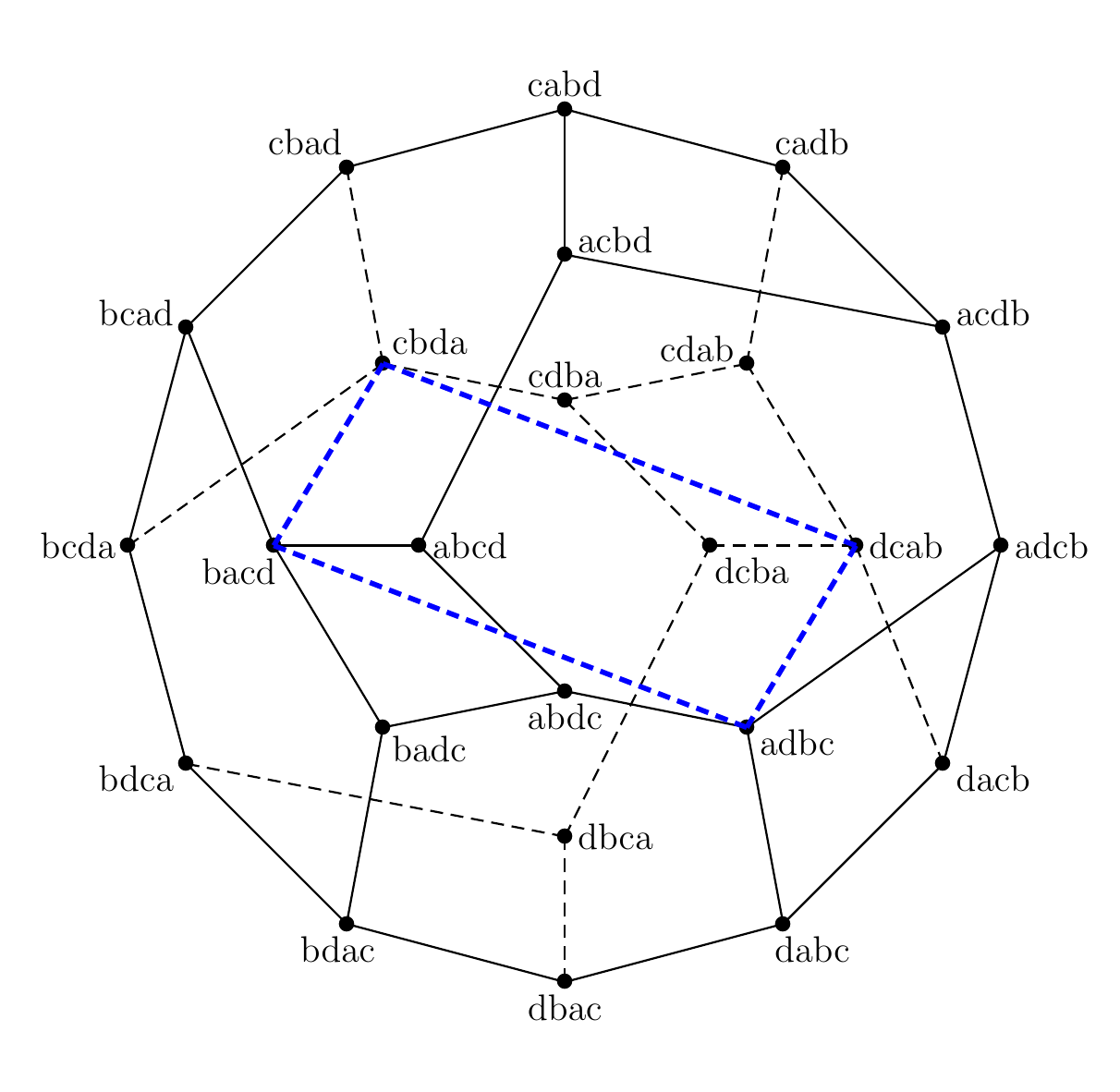}
\caption{A maximal Condorcet domain with  four elements}
\label{Fig9}
\end{figure}

%\begin{center} 
%\includegraphics[width=9cm]{Figure_9}
%\vspace*{1mm} \\
%{\em Figure 9: A maximal Condorcet domain with four elements}
%\vspace{3mm}
%\end{center}

%%%%%%%%%%%%%%%%%%%%%%%%%%%%%%%%%%%%%%%%%%%%%%%
%%%%%%%%%%%%%%%%%%%%%%%%%%%%%%%%%%%%%%%%%%%%%%%
%%% SECTION 6
%%%%%%%%%%%%%%%%%%%%%%%%%%%%%%%%%%%%%%%%%%%%%%%
%%%%%%%%%%%%%%%%%%%%%%%%%%%%%%%%%%%%%%%%%%%%%%%

\section{Arrovian Aggregation and Strategy-Proof Social Choice on
Median Preference Domains}
\setcounter{equation}{0}

Closed Condorcet domains not only preclude intransitivities in pairwise majority voting,
they are also endowed with a large
class of further aggregation rules satisfying Arrow's independence condition.
This follows from the analysis of \cite{NehringPuppe2007b,NehringP2010}. 
Indeed, their
main result entails a characterization of all Arrovian aggregators on such
domains under an additional monotonicity
condition. In the first subsection below, we apply their result to describe the class of
all monotone Arrovian aggregators on closed Condorcet domains.
The monotonicity condition plays then a crucial role in the construction
of strategy-proof social choice functions on such domains in the second
subsection.

\subsection{Characterization of all Arrovian Aggregators}
Let $N=\{ 1,2,\ldots ,n\}$ be the set of voters.
An {\em aggregator} on a domain $\cd\subseteq\cR (X)$ is a mapping
$f\colon \cd^n\rightarrow \cd$ that assigns an order in $\cd$ to each
profile of individual orders in $\cd$. The following
conditions on aggregators have been extensively studied in the literature.\vspace{1ex}\\
{\bf Full Range.} For all $R\in\cd$, there exist $R_1,\ldots,R_n$ such that
$f(R_1,\ldots,R_n)=R$.\vspace{1mm}\\
{\bf Unanimity.} For all $R\in \cd$, $f(R,\ldots,R)=R$.\vspace{1mm}\\
{\bf Independence.} For all $R_1,\ldots,R_n, R'_1,\ldots,R'_n\in \cd$ and all pairs
of distinct alternatives $x,y\in X$,
if $xRy$ and, for all $i\in N$,
[$xR_iy\Leftrightarrow xR'_iy$], then $xR'y$,
where $R=f(R_1,\ldots,R_n)$ and $R'=f(R'_1,\ldots,R'_n)$.\vspace{1ex}

An aggregator is called {\dfn Arrovian} if it satisfies unanimity and
independence. In what follows, we will be concerned with Arrovian
aggregators that satisfy in addition the following monotonicity
condition.\vspace{1ex}\\
{\bf Monotonicity} For all $R_1,\ldots,R_n, R'_i\in \cd$ and all pairs
of distinct alternatives $x,y$,
if $xR'y$ and $x R_i y$, then
$xRy$, where $R' = f(R_1,\ldots,R'_i,\ldots,R_n)$ and $R=f(R_1,\ldots,R_i,\ldots,R_n)$.
\vspace{1ex}

The monotonicity condition can be rephrased as follows. An aggregator
$f$ is {\dfn monotone} if and only if, for all
$R_i,R'_i\in \cd$ and all $R_{-i}\in \cd^{n-1}$,
\begin{equation}
\label{monotone}
f(R_i,R_{-i})\in [ R_i,f(R'_i,R_{-i})] \ .
\end{equation}
This has the following interpretation. Consider any pair of distinct alternatives
$a$ and $b$, and suppose that $aR_i b$ according to agent $i$'s true order $R_i$.
Then, if agent $i$ can force the social order to rank $a$
above $b$ by submitting some order $R'_i$, the social order would also
rank $a$ above $b$ if agent $i$ submitted his true preference $R_i$. In other words,
no agent can benefit  in a pairwise comparison from any
misrepresentation. The monotonicity condition thus
has a clear
`non-manipulability' flavor which we will further exploit. 

The conjunction of independence and
monotonicity is equivalent to the following single condition.
\vspace{1ex}\\
{\bf Monotone Independence} For all $R_1,\ldots,R_n$ and for all pairs $(x,y)$ of distinct
alternatives $x,y\in X$, if $xRy$, where $R=f(R_1,\ldots,R_n)$ and
[$a R_i b\Rightarrow a R'_ib$] for all $i\in N$,
then $aR'b$, where $R'=f(R'_1,\ldots,R'_n)$.\vspace{1ex}

Note also that under monotonicity, unanimity can be deduced
from the full range condition, i.e., from the assumption that the aggregator is
onto.
\medskip

For every pair $(x,y)\in X\times X$ of distinct alternatives, let $\cw_{xy}$ be a non-empty
collection of non-empty subsets $W\subseteq N$ of voters
(the {\em winning coalitions for $x$ against $y$} ) satisfying
\begin{equation}
\left[ W\in\cw_{xy}\mbox{ and }W'\supseteq W\right] \ \Rightarrow \
W'\in\cw_{xy} .
\end{equation}
\begin{definition}
A collection
$\cw =\left\{ \cw_{ab}  \mid  (a,b)\in X\times X, a\neq b \, \right\}$ is called a
{\dfn structure of winning coalitions} if, for all distinct pairs of alternatives
$(x,y)\in X\times X$,
\begin{equation}
\label{proper}
W\in\cw_{xy}\ \Leftrightarrow \ W^c\not\in\cw_{yx},
\end{equation}
where $W^c=N\setminus W$ denotes the complement of the coalition $W$.
\end{definition}

\noindent
{\bf Examples.} There are numerous examples of structures of winning
coalitions. The simplest are the {\em dictatorships} which arise whenever
there exists an individual $i$ such that, for all distinct pairs $(x,y)$, a coalition
$W$ belongs to $\cw_{xy}$ if and only if $i\in W$. More generally, an {\em oligarchy}
is characterized by the existence of a non-empty set $M$ of voters such
that, for all distinct pairs $(x,y)$, either (i) $W\in\cw_{xy}\Leftrightarrow M\subseteq W$,
or (ii) $W\in\cw_{yx}\Leftrightarrow M\subseteq W$. Note that, by \eqref{proper}, if (i)
holds, then $\{i\}\in\cw_{yx}$ for all $i\in M$; similarly, if (ii) holds,
then $\{i\}\in\cw_{xy}$ for all $i\in M$.

A structure of winning coalitions is {\em anonymous} if, for any fixed pair $(x,y)\in X\times X$,
membership of a coalition $W$ to $\cw_{xy}$  depends only on the number
of individuals in $W$. A special case is {\em pairwise majority voting} which requires
that for all pairs $(x,y)\in X\times X$ of distinct alternatives
\[
W\in\cw_{xy} \ \Leftrightarrow \  |W| > n/2.  
\]
Note that this is jointly satisfiable with \eqref{proper} only if $n$ is odd.
\medskip

\noindent
Given a domain $\cd\in\cR(X)$, a structure of winning coalitions is said to be
{\dfn order preserving on $\cd$} if,
for every pair of distinct alternatives $(x,y)\in X\times X$ and every pair of distinct
alternatives $(z,w)\in X\times X$,
\begin{equation}
\cvd_{xy}\subseteq \cvd_{zw} \ \Rightarrow \
\cw_{xy} \subseteq \cw_{zw}.
\end{equation}
Observe that pairwise majority voting always defines an order preserving structure
of winning coalitions (since $\cw_{xy} = \cw_{zw}$ for all distinct $x,y$
and $z,w$ in the case of pairwise majority voting).

The following result is due to \cite{NehringPuppe2007b}.
\begin{proposition}
\label{charctMAR} Let $\cd$ be a closed Condorcet domain and let $\cw$ be any
order preserving structure of winning coalitions on $\cd$. For all preference
profiles $(R_1,\ldots,R_n)\in\cd^{n}$ there exists a unique order $R^*\in \cd$ such that,
for every pair of distinct alternatives ${(x,y)\in X\times X}$
\begin{equation}
\label{defR*}
xR^* y \ \Leftrightarrow \ \{ i\in N \mid xR_i y\}\in\cw_{xy}\,.
\end{equation}
The aggregator defined by \eqref{defR*}  is a monotone Arrovian aggregator. Conversely,
every monotone Arrovian aggregator on $\cd$ takes this form for some order preserving
structure of winning coalitions $\cw$.
\end{proposition}
{\bf Proof.} The proof can be deduced from \cite[Prop.~3.4]{NehringPuppe2007b}. 
For the sake of the paper being self-contained, we reproduce the argument here, again
emphasizing the role of the Helly property on median domains as in the
proof of Theorem 2 above.

For any profile $(R_1,\ldots,R_n)\in\cd^{n}$, define a binary relation $R^*$ by \eqref{defR*}.
First, we show that $R^*\in\cd$. Consider distinct $x,y\in X$ and
distinct $z,w\in X$ such that $xR^* y$ and $zR^*w$; we will show
that $\cvdxy\cap\cvd_{zw} \ne \emptyset$. By contradiction, suppose that this
does not hold, then we would obtain $\cvdxy\subseteq\cvd_{wz}$, and in particular
also $\{ i\in N \mid xR_i y\}\subseteq \{ i\in N \mid wR_i z\}$. At the same time,
we would have $\cw_{xy}\subseteq\cw_{wz}$ because the structure of winning
coalitions is order preserving. Thus, $\{ i\in N \mid xR_i y\}\in\cw_{wz}$,
hence also $\{ i\in N \mid wR_i z\}\in\cw_{wz}$. But the latter contradicts
$zR^*w$ using \eqref{proper}. Thus, the collection of convex sets $\{ \cvdxy \mid xR^*y\}$
has pairwise non-empty intersections. By the Helly property for convex sets in any
Condorcet domain (Proposition~\ref{helly-prop}),
we obtain
\[
\bigcap_{\{ x R^* y\} }\cvdxy \neq\emptyset,
\]
which means that $R^*$ is an element of $\cd$.
Thus, \eqref{defR*}  indeed defines a mapping from $\cd^n$ to $\cd$. It is easily
verified that it is onto, monotone and independent, i.e., an Arrovian aggregator.

Conversely, let $f$ be a monotone Arrovian aggregator. By the monotonicity and
independence conditions, it is easily seen that $f$ can be defined in terms of a structure
of winning coalitions as in \eqref{defR*}. We thus have only to verify that the corresponding
structure of winning coalitions must be order preserving. Thus, assume that, for
distinct $(x,y)\in X\times X$ and distinct $(z,w)\in X\times X$, we have
$\cvd_{xy}\subseteq \cvd_{zw}$, and $W\in\cw_{xy}$. Then, if the
profile $(R_1,\ldots,R_n)$ is such that all voters in $W$ prefer $x$ to $y$, we must
have $xRy$ where $R=f(R_1,\ldots,R_n)$. By assumption, all orders in $\cd$ which
rank $x$ above $y$ must also rank $z$ above $w$, hence $zRw$. By independence,
this holds for {\em all} profiles in which the agents in $W$ rank $x$ above $y$, hence
also $z$ above $w$, i.e., the agents in $W$ are also winning for $z$ versus $w$.
%\vspace{-1mm} \\  
\hspace*{1mm}\hfill  $\Box$

%%%%%%%%%%%%%%%%%%%%%%%%%%%%%%%%%%%%%
%%%%%% SECTION 6.2
%%%%%%%%%%%%%%%%%%%%%%%%%%%%%%%%%%%%%

\subsection{Strategy-Proof Social Choice}
It is well-known that on domains on which pairwise majority voting with an odd number
of voters is transitive,
choosing the Condorcet winner yields a strategy-proof social choice function
(see, e.g., Lemma 10.3 in \cite{Moulin1988}). In this final
subsection, we use Proposition~\ref{charctMAR} and property \eqref{monotone}  which is
entailed by monotonicity to construct a rich class of further
strategy-proof social choice functions on any closed Condorcet domain.

A social choice function $F$ that maps every profile $(R_1, \ldots ,R_n)\in\cd^n$
to an element $F (R_1, \ldots ,R_n)\in X$ is {\dfn strategy-proof} if, for all $i\in N$,
all $R_i,R'_i\in \cd$ and all $R_{-i}\in \cd^{n-1}$,
\[
[F(R_i,R_{-i})]  R_i  [F(R'_i,R_{-i})],
\]
i.e., if no voter can benefit by misrepresenting her true preferences.

For each order $R\in\cR (X)$ denote by $\tau (R)\in X$ the top-ranked
element of~$R$. Let $\cd\subseteq\cR(X)$ be any closed Condorcet domain, and
consider any order preserving structure of winning coalitions $\cw$.
For every profile $(R_1, \ldots ,R_n)\in\cd^n$ let $R_{\cw}\in\cd$ be the unique
order satisfying  \eqref{defR*}  for all distinct $x,y\in X$. Define a social choice
function $F_{\cw}\colon \cd^n\rightarrow X$ by
\begin{equation}
\label{defFW}
F_{\cw}(R_1,...,R_n) = \tau (R_{\cw}).
\end{equation}
\begin{theorem}
\label{strat-proof} Let $\cd\subseteq\cR(X)$ be any closed Condorcet domain. For every
order preserving structure of winning of coalitions $\cw$, the social choice function
$F_{\cw}$ defined by \eqref{defFW}  is strategy-proof.
\end{theorem}

{\bf Proof.} By Proposition~\ref{charctMAR}, the aggregator
$f_{\cw}\colon \cd^n\rightarrow\cd$ that maps
any profile $(R_1,...,R_n)$ to the social order $R_{\cw}$ according to \eqref{defR*} 
is a monotone Arrovian aggregator; in particular, it satisfies \eqref{monotone}.
In other words, if we denote $R_{\cw} = f_{\cw} (R_i,R_{-i})$
and $R'_{\cw} = f_{\cw} (R'_i,R_{-i})$, we have for all $R_1,...,R_n$, all $i$,
all $R'_i$, and all distinct pairs of alternatives $(x,y)\in X\times X$,
\begin{equation}
\label{impFW}
\left[ xR_i y\mbox{ and }xR'_{\cw}y \right] \ \Rightarrow \ xR_{\cw}y
\end{equation}
This implies at once the strategy-proofness of $F_{\cw}$, by contraposition.
Indeed, suppose that agent $i$ could benefit by misreporting $R'_i$, i.e., suppose
that $\tau (R'_{\cw})R_i\tau (R_{\cw})$, where $R_i$ is agent $i$'s true
preference order. Then, since
$\tau (R'_{\cw})$ is the top element
of the order $R'_{\cw}$ we obtain from \eqref{impFW} ,
$\tau (R'_{\cw})R_{\cw}\tau (R_{\cw})$. Since
$\tau (R_{\cw})$ is the top element of $R_{\cw}$
this implies $\tau (R'_{\cw}) = \tau (R_{\cw})$, i.e., the
misrepresentation does not change the chosen alternative.
%\vspace{-1mm} \\  
\hspace*{1mm}\hfill \vspace{2ex} $\Box$

In case of the classical single-crossing property, the anonymous
social choice functions defined by \eqref{defFW}  are exactly the ones identified by
\cite{Sapor2009}. Thus, in this special case, the class of  anonymous social
choice functions
considered in Theorem \ref{strat-proof} exhausts all  anonymous
strategy-proof social choice
functions. It is an open {\color{black} and interesting} question whether this holds more generally on all
closed Condorcet domains  (and whether the anonymity assumption is really
necessary for this conclusion).

%%%%%%%%%%%%%%%%%%%%%%%%%%%%%%%%%%%%%%%%%%%%%%%%
%%%%%%%%%%%%%%%%%%%%%%%%%%%%%%%%%%%%%%%%%%%%%%%%
%%%%%%%%%%%%%%%%%%%%%%%%%%%%%%%%%%%%%%%%%%%%%%%%
%%%%%%%%%%%%%%%%%%%%%%%%%%%%%%%%%%%%%%%%%%%%%%%%
%%%%%%%%%%%%%%%%%%%%%%%%%%%%%%%%%%%%%%%%%%%%%%%%
%%%%%%%%  APPENDIX
%%%%%%%%%%%%%%%%%%%%%%%%%%%%%%%%%%%%%%%%%%%%%%%%
%%%%%%%%%%%%%%%%%%%%%%%%%%%%%%%%%%%%%%%%%%%%%%%%
%%%%%%%%%%%%%%%%%%%%%%%%%%%%%%%%%%%%%%%%%%%%%%%%
%%%%%%%%%%%%%%%%%%%%%%%%%%%%%%%%%%%%%%%%%%%%%%%%
%%%%%%%%%%%%%%%%%%%%%%%%%%%%%%%%%%%%%%%%%%%%%%%%

%\end{document}

\vspace{1cm}
%\newpage
%
%
\renewcommand{\thesection}{A}
\begin{center}
{\Large {\bf Appendix: Remaining Proofs}}
\end{center}
\setcounter{theorem}{0}
\setcounter{equation}{0}
\setcounter{proposition}{0}
Lemma \ref{equivalence} is proven using a result by \cite{BandeltChepoi1996}. The
statement of this result
requires some additional definitions, in particular the notion of a `geometric interval operator,'
as follows. An {\dfn interval operator}
on a (finite) set $V$ is a mapping that assigns to each pair $(v,w)\in V\times V$ a non-empty subset
$[v,w]\subseteq V$, the {\em interval spanned by $v$ and $w$}, such that, for all $v,w\in V$,
$v\in [v,w]$ and $[v,w] = [w,v]$. An interval operator is called {\dfn geometric} if it satisfies in addition
the following properties. For all $t,u,v,w\in V$,
\begin{eqnarray}
 [v,v] & = & \{ v\} , \\
 u\in [v,w] & \Rightarrow & [v,u]\subseteq [v,w] , \\
 t,u\in [v,w] \ \& \ t\in [v,u] & \Rightarrow & u \in [t,w] .
\end{eqnarray}

A pair $v,w$ is called an {\dfn edge} if $v\neq w$ and $[v,w]=\{ v,w\}$. These edges form
a graph $\Gamma$ on the vertex set $V$.

\begin{lemma}\label{max_chain} Consider a
geometric interval operator on $V$, and let $u\in [v,w]$. Then, there
exist pairwise distinct $t_1, \ldots , t_m \in [v,w]$ such that $t_1 = v$, $t_m =w$, $t_k = u$ for some
$k\in \{ 1, \ldots , m\}$, and such that
\[
[t_1, t_2]\subset [t_1,t_3]\subset \ldots \subset [t_1,t_m]
\]
forms a maximal chain. In particular, the graph induced by a geometric interval operator
is connected.
\end{lemma}
\noindent {\bf Proof.} The existence of a maximal chain of the required form follows at once from
condition (A.2). The pairs $t_kt_{k+1}$ must form an edge for $k=1, \ldots m-1$ by
maximality of the chain. Thus, any two vertices are connected by a path.
%\vspace{-1mm} \\  
\hspace*{1mm} %\vspace{2ex}
\hfill $\Box$

\medskip

\noindent An interval operator is called {\dfn graphic} if, for all $u,v,w\in V$, $u\in [v,w]$
if and only if $u$ is geodesically between $v$ and $w$ in the induced graph $\Gamma$; note
that this is exactly condition (\ref{equivalence_equation}) in Lemma \ref{equivalence} above. An
interval operator is said to satisfy the
{\dfn triangle condition} if, for all triples $u,v,w\in V$ such that
\begin{equation}
[u,v] \cap [v,w] = \{ v\} \ \mbox{ and } \ [v,w] \cap [w,u] = \{ w\} \
 \mbox{ and } \ [w,u] \cap [u,v] = \{ u\} ,
\end{equation}
all three intervals are edges whenever one of them is. Observe that (A.4)
can be satisfied only if either all three elements $u,v,w$ coincide, or are
pairwise distinct. The following result is due to
\cite[Th.1]{BandeltChepoi1996}.

\begin{proposition}\label{triangle} Any geometric
interval operator satisfying the triangle condition is graphic.
\end{proposition}

We want to apply Proposition \ref{triangle} to prove Lemma \ref{equivalence}. In order
to do so, we first verify the geometricity of the natural interval operator induced by every
subdomain of orders. 
 
\begin{lemma} \label{geometric} For any domain $\cd\subseteq \cR (X)$, the interval operator
that assigns to every pair $R,R'\in\cd$ the interval $[R,R']\cap\cd$ is geometric. 
\end{lemma} 
\noindent {\bf Proof.} Properties (A.1) and (A.2) are easily verified. To verify the so-called `inversion
law' (A.3), consider $T,U,V,W\in\cd$ as required in the antecedent of (A.3). Let $x,y\in X$ be such
that $xTy$ and $xWy$. We have to show that then $xUy$. Since $x,y$ were arbitrarily chosen,
this would imply $U\in[T,W]$, as desired. Assume, by contradiction, that $yUx$; then, we must
have $xVy$ since by assumption $T\in [V,U]$. But this contradicts the assumption that
$U\in [V,W]$.%\vspace{-1mm} \\  
\hspace*{1mm}% \vspace{2ex}
\hfill $\Box$

\medskip

The following proof of Lemma \ref{equivalence} shows that the triangle condition is
a powerful sufficient condition for an interval operator to be graphic, since it is indeed
satisfied in all three cases considered in Lemma \ref{equivalence} (sometimes vacuously).
\vspace{2ex}

\noindent
{\bf Proof of Lemma \ref{equivalence}. (i)} In case of a median domain, the triangle condition is vacuously
satisfied, since there can obviously be no triples of pairwise distinct elements satisfying (A.4)
by the median property. By Proposition \ref{triangle} the equivalence (\ref{equivalence_equation})
is satisfied for any median domain.

{\bf (ii)} Next consider any connected domain $\cd\subseteq\cR (X)$. There can exist
triples satisfying (A.4), but we shall show that in this case
none of the three intervals can form an edge,
hence the triangle condition is again satisfied. Thus, suppose that the three pairwise
distinct orders
$U,V,W\in\cd$ satisfy (A.4) and, by contradiction, that one of the three intervals is an edge,
say $[U,V]=\{ U,V\}$. Since $\cd$ is connected, there exists $x,y\in X$ such that
$U$ and $V$ differ only in the ranking of $x$ versus $y$, say $xUy$ and $yVx$,
while $U$ and $V$ agree in the ranking of all other pairs of alternatives. There are
two possibilities: either $xWy$ or $yWx$. In the first case, we have $U\in [V,W]$ and 
hence $([U,V] \cap [V,W]) \supseteq \{ U,V\}$; in the second case, $V\in [W,U]$ and hence
$([W,U] \cap [U,V]) \supseteq \{ U,V\}$. In both cases, we thus obtain a contradiction to
assumption (A.4). By Proposition \ref{triangle} the equivalence (\ref{equivalence_equation})
is satisfied for $\cd$.

{\bf (iii)} Finally, assume that $\cd$ is such that $\Gamma_{\cd}$ is acyclic, i.e., a tree.
As in part (i), we show that there cannot exist triples satisfying (A.4) hence again the triangle
condition is satisfied vacuously.\footnote{Note, however, that we cannot use part (i) directly since
we do not know yet whether $\cd$ is a median domain.} Assume, by way of contradiction,
that the pairwise distinct orders $U,V,W\in\cd$ satisfy (A.4). By Lemma \ref{max_chain},
there exists a path $\pi_{UV}$ in $\Gamma_{\cd}$ connecting $U$ and $V$ that stays
entirely in $[U,V]$; in particular, $\pi_{UV}$ does not contain $W$. Similarly,
there exists a path $\pi_{VW}$ in $\Gamma_{\cd}$ connecting $V$ and $W$
and staying entirely in $[V,W]$, and a path $\pi_{WU}$ connecting $W$ and $U$
and staying entirely in $[W,U]$.
But then the union $\pi_{UV}\cup \pi_{VW}\cup \pi_{WU}$ forms a cycle, which contradicts
the assumed acyclicity of $\Gamma_{\cd}$. Thus, again,
by Proposition \ref{triangle} the equivalence (\ref{equivalence_equation})
is satisfied for the domain $\cd$.
%\vspace{-1mm} \\  
\hspace*{1mm}\hfill  $\Box$

\medskip

Note that there do exist connected domains $\cd$ and pairwise distinct orders
$U,V,W\in\cd$ satisfying condition (A.4). Examples are all triples of orders on a
common $6$-cycle with pairwise distance of two such as, for instance,
$abc$, $bca$, $cab$ in Figure~2, or $abcd$, $bcad$, $cabd$ in Figure~5, etc.

\vspace{1cm}
%\newpage

\bibliographystyle{plainnat}
\bibliography{cps}

\end{document}